\documentstyle[12pt]{article}

\begin{document}
\medskip

\rightline {math.QA/9805048 \hskip 10 mm}
\vskip 50 pt

\noindent {\bf REPRESENTATIONS OF THE CYCLICALLY SYMMETRIC}

\noindent{\bf $q$-DEFORMED ALGEBRA ${\bf so}_{\bf q}({\bf 3})$}

\vskip 15 pt

{\sl M. Havl\'{i}\v cek}

{\it Department of Mathematics, FNSPE, Czech Technical University}

{\it CZ-120 00, Prague 2, Czech Republic}
\medskip

{\sl A. U. Klimyk}

{\it Institute for Theoretical Physics, Kiev 252143, Ukraine}
\medskip

{\sl S. Po\v sta}

{\it Department of Mathematics, FNSPE, Czech Technical University}

{\it CZ-120 00, Prague 2, Czech Republic}

\vskip 30 pt

\begin{abstract}
An algebra homomorphism $\psi$ from the nonstandard $q$-deformed (cyclically
symmetric) algebra $U_q({\rm so}_3)$ to the extension
${\hat U}_q({\rm sl}_2)$ of the Hopf algebra $U_q({\rm sl}_2)$ is
constructed. Not all irreducible representations of
$U_q({\rm sl}_2)$ can be extended to representations of
${\hat U}_q({\rm sl}_2)$. Composing the homomorphism $\psi$ with irreducible
representations of ${\hat U}_q({\rm sl}_2)$ we obtain representations of
$U_q({\rm so}_3)$. Not all of these representations of $U_q({\rm so}_3)$
are irreducible. Reducible representations of $U_q({\rm so}_3)$
are decomposed into irreducible components. In this way we obtain all
irreducible representations of $U_q({\rm so}_3)$ when $q$ is not a root of
unity. A part of these
representations turns into irreducible representations of the Lie algebra
so$_3$ when $q\to 1$. Representations of the other part have no
classical analogue. Using the homomorphism $\psi$ it is shown how
to construct tensor products of finite dimensional representations of
$U_q({\rm so}_3)$. Irreducible representations of $U_q({\rm so}_3)$
when $q$ is a root of unity are constructed. Part of them are
obtained from irreducible representations of ${\hat U}_q({\rm sl}_2)$
by means of the homomorphism $\psi$.

\end{abstract}

\newpage

\noindent
{\sf I. INTRODUCTION}
\bigskip

It is well-known that the Lie algebras ${\rm sl}_2$ and ${\rm so}_3$ of the
Lie groups $SL(2,{\bf C})$ and $SO(3)$, respectively, are isomorphic. But
these algebras differ from each other if we consider their embedding to
the wider Lie algebra ${\rm sl}_3$. There is no automorphism of
${\rm sl}_3$ which transfers the embedding ${\rm sl}_2 \subset {\rm sl}_3$
to the embedding ${\rm so}_3 \subset {\rm sl}_3$. Note that the
embedding ${\rm so}_3 \subset {\rm sl}_3$ is of great
importance for nuclear physics:
it is used in spectroscopy.

The definition of the $q$-analogue of the universal enveloping algebra
$U({\rm sl}_2)$ is well-known. It is the quantum algebra
$U_q({\rm sl}_2)$ which is a Hopf algebra. If we wish to have a $q$-analogue
of the universal enveloping algebra ${\rm so}_3$ such that at $q\to 1$
we obtain the classical embedding ${\rm so}_3 \subset {\rm sl}_3$,
then the algebra ${\rm sl}_2$ is not appropriate for this role. By other
words, an algebra $U_q({\rm so}_3)$ must differ from $U_q({\rm sl}_2)$.
This algebra $U_q({\rm so}_3)$ is well. It is the associative
algebra generated by three elements $I_1$, $I_2$ and $I_3$ satisfying
the relations
$$
q^{1/2}I_1I_2-q^{-1/2}I_2I_1=I_3,  \eqno (1)$$
$$
q^{1/2}I_2I_3-q^{-1/2}I_3I_2=I_1, \eqno (2)$$
$$
q^{1/2}I_3I_1-q^{-1/2}I_1I_3=I_2. \eqno (3)
$$
Such (and more general) deformation of the commutator $[I_i,I_j]=
I_iI_j-I_jI_i$ was defined at 1967 by R. Santilli in the paper [1] (see
also [2] and [3]) under studying a generalization of the Lie theory.
Afterwards (in 1990), the algebra $U_q({\rm so}_3)$ with commutation
relations (1)--(3) was determined by D. Fairlie [4]. An algebra which
can be reduced to $U_q({\rm so}_3)$ was defined in 1986 by M.
Odesski [5].

Fairlie [4] gave finite dimensional irreducible representations of the
algebra $U_q({\rm so}_3)$ which at $q\to 1$ give the well-known finite
dimensional irreducible representations of the Lie algebra so$_3$. These
representations are given by integral or half-integral non-negative numbers.
Odesski [5] also gave some classes of irreducible representations.

It was shown (see [5--7]) that the algebra $U_q({\rm so}_3)$
has irreducible finite dimensional representations which have no classical
analogue (that is, which do not admit the limit $q\to 1$). It was not clear
why such strange representations of the algebra $U_q({\rm so}_3)$ appear.
What is their nature? The answer to this question is one of the aims of
this paper.

We construct a homomorphism from $U_q({\rm so}_3)$ to the algebra
${\hat U}_q({\rm sl}_2)$ which is an extension of the well-known quantum
algebra $U_q({\rm sl}_2)$ (note that there is no homomorphism from
$U_q({\rm so}_3)$ to $U_q({\rm sl}_2)$).
Irreducible finite dimensional representations of
$U_q({\rm sl}_2)$ (but not all) can be extended to finite dimensional
representations of the algebra ${\hat U}_q({\rm sl}_2)$. Composing
a homomorphism $U_q({\rm so}_3)\to {\hat U}_q({\rm sl}_2)$ with
these representations of ${\hat U}_q({\rm sl}_2)$, we obtain
representations of the algebra $U_q({\rm so}_3)$. But some of irreducible
representations of ${\hat U}_q({\rm sl}_2)$ lead to reducible representations
of the algebra $U_q({\rm so}_3)$. Decomposing these reducible
representations of $U_q({\rm so}_3)$ we obtain irreducible representations of
this algebra which have no analogue for the Lie algebra so$_3$.
If $q$ is not a root of unity, then in this way we obtain all finite
dimensional irreducible representations of $U_q({\rm so}_3)$. But there are
infinite dimensional irreducible representations of $U_q({\rm so}_3)$
which cannot be obtained in this way.

Existence of the homomorphism $U_q({\rm so}_3)\to {\hat U}_q({\rm sl}_2)$
allows us to define tensor products of representations of the algebra
$U_q({\rm so}_3)$ which is not a Hopf algebra.

Using the homomorphism $U_q({\rm so}_3)\to {\hat U}_q({\rm sl}_2)$
and irreducible representations of ${\hat U}_q({\rm sl}_2)$ we obtain
representations of $U_q({\rm so}_3)$ when $q$ is a root of unity.
Taking irreducible representations of $U_q({\rm so}_3)$ obtained in this way
and decomposing reducible representations, we obtain several series of
irreducible representations of $U_q({\rm so}_3)$. In addition, we construct
irreducible representations of $U_q({\rm so}_3)$ which cannot be derived
from ${\hat U}_q({\rm sl}_2)$.

When $q$ is not a root of unity, then each irreducible (finite or
infinite dimensional) representation of $U_q({\rm so}_3)$ is equivalent
to one of the representations constructed below. (We do not give a proof
of this assertion in this paper because it would take much place; this
proof will be given in a separate paper.) We think that in this paper we
constructed also all irreducible representations of $U_q({\rm so}_3)$
when $q$ is a root of unity. But in this case we have no proof of this
assertion.
The reason of this is that in this case there are many classes of irreducible
representations and a proof of completeness of irreducible representations
becomes very tedious.

Let us remark that in [5] there were constructed irreducible finite
dimensional representations of $U_q({\rm so}_3)$ when $q$ is not a root
of unity and a part of irreducible infinite dimensional representations.
In [6] and [7], there were constructed irreducible representations
of $U_q({\rm so}_3)$ which satisfy the conditions of $*$-representations
(that is, such that $T(I^*_j)=-T(I_j)$, $j=1,2$). These $*$-representations
are a part of irreducible representations of $U_q({\rm so}_3)$
constructed in this paper. We started to study irreducible representations of
$U_q({\rm so}_3)$ for $q$ a root of unity in [8], where a part of
irreducible representations for this case were constructed.
Note that in [5--8] there are no relations of
representations of $U_q({\rm so}_3)$ to representations of
${\hat U}_q({\rm sl}_2)$. This relation makes representations of
$U_q({\rm so}_3)$ clear and understandable.

We suppose that in Sections II and III $q$ is any complex number different
from --1. In Sections IV--VII, $q$ is not a root of unity. In
Sections VIII--X, $q$ is a root of unity.
\bigskip

\noindent
{\sf II. THE ALGEBRAS $U_q({\rm so}_3)$ AND ${\hat U}_q({\rm sl}_2)$}
\medskip

The algebra $U_q({\rm so}_3)$ is obtained by a $q$-deformation of the
standard commutation relations
$$
[I_1,I_2]=I_3,\ \ \ [I_2,I_3]=I_1,\ \ \ [I_3,I_1]=I_2
$$
of the Lie algebra so$_3$. So, $U_q({\rm so}_3)$ is defined as the
complex associative algebra with unit element generated by the elements
$I_1$, $I_2$, $I_3$ satisfying the defining relations
$$
[I_1,I_2]_q:= q^{1/2}I_1I_2-q^{-1/2}I_2I_1=I_3,  \eqno (4)$$
$$
[I_2,I_3]_q:= q^{1/2}I_2I_3-q^{-1/2}I_3I_2=I_1, \eqno (5)$$
$$
[I_3,I_1]_q :=q^{1/2}I_3I_1-q^{-1/2}I_1I_3=I_2. \eqno (6)
$$
Unfortunately, a Hopf algebra structure is not known on $U_q({\rm so}_3)$.
However, it can be embedded into the Hopf algebra $U_q({\rm sl}_3)$ as a
Hopf coideal (see [9]). This embedding is very important for the possible
application in spectroscopy.

It follows from the relations (4)--(6) that for the algebra $U_q({\rm so}_3)$
the Poincar\'e--Birkhoff--Witt theorem is true and this theorem can be
formulated as: {\it The elements $I_1^kI_2^mI_3^n$, $k,m,n=0,1,2,\cdots $,
form a basis of the linear space $U_q({\rm so}_3)$.} Indeed, by using
the relations (4)--(6) any product $I_{j_1}I_{j_2}\cdots I_{j_s}$,
$j_1,j_2,\cdots ,j_s=1,2,3$, can be reduced to a sum of the elements
$I_1^kI_2^mI_3^n$ with complex coefficients.

Note that by (4) the element $I_3$ is not independent: it is determined
by the elements $I_1$ and $I_2$. Thus, the algebra $U_q({\rm so}_3)$
is generated by $I_1$ and $I_2$, but now instead of quadratic relations
(4)--(6) we must take the relations
$$
I_1I_2^2-(q+q^{-1})I_2I_1I_2+I_2^2I_1=-I_1, \eqno (7)$$
$$
I_2I_1^2-(q+q^{-1})I_1I_2I_1+I_1^2I_2=-I_2, \eqno (8)
$$
which are obtained if we substitute the expression (4) for $I_3$ into
(5) and (6). The equation $I_3=q^{1/2}I_1I_2-q^{-1/2}I_2I_1$ and the
relations (7) and (8) restore the relations (4)--(6).

Remark that the definition of $U_q({\rm so}_3)$ by means of relations
(7) and (8) was used in [9] for the embedding of $U_q({\rm so}_3)$ to
$U_q({\rm sl}_3)$. The relations (7) and (8) differ from Serre's relations
in the definition of quantum algebras by V. Drinfeld and M. Jimbo by
appearance of non-vanishing right hand sides.

The algebra $U_q({\rm so}_3)$ is closely related to (but not coincides
with) the quantum algebra $U_q({\rm sl}_2)$. The last algebra is generated
by the elements $q^H$, $q^{-H}$, $E$, $F$ satisfying the relations
$$
q^Hq^{-H}=q^{-H}q^H=1,\ \ \ q^HEq^{-H}=qE,\ \ \ q^HFq^{-H}=q^{-1}F, \eqno (9)
$$
$$
[E,F]:= EF-FE=\frac {q^{2H}-q^{-2H}}{q-q^{-1}} . \eqno (10)
$$
Note that $U_q({\rm sl}_2)$ is the associative algebra equipped with a Hopf
algebra structure (a comultiplication, a counit and an antipode). In
particular, the comultiplication $\Delta$ is determined by the formulas
$$
\Delta (q^{\pm H})=q^{\pm H}\otimes q^{\pm H},\ \ \
\Delta (E)=E\otimes q^{H}+q^{-H}\otimes E, $$
$$
\Delta (F)=F\otimes q^{H}+q^{-H}\otimes F.
$$

In order to relate the algebras $U_q({\rm so}_3)$ and $U_q({\rm sl}_2)$
we need to extend $U_q({\rm sl}_2)$ by the elements
$(q^kq^H+q^{-k}q^{-H})^{-1}$ in the sense of [10].
We denote by ${\hat U}_q({\rm sl}_2)$ the associative algebra with unit
element generated by the elements
$$
q^H,\ \  \ q^{-H},\ \ \ E,\ \ \ F,\ \ \ (q^kq^H+q^{-k}q^{-H})^{-1},
\ \ k\in {\bf Z},
$$
satisfying the defining relations (9) and (10)
of the algebra $U_q({\rm sl}_2)$ and
the following natural relations:
$$
(q^kq^H+q^{-k}q^{-H})^{-1}(q^kq^H+q^{-k}q^{-H})=
(q^kq^H+q^{-k}q^{-H})(q^kq^H+q^{-k}q^{-H})^{-1}=1, \eqno (11) $$
$$
q^{\pm H}(q^kq^H+q^{-k}q^{-H})^{-1}=(q^kq^H+q^{-k}q^{-H})^{-1}q^{\pm H} ,
\eqno (12) $$
$$
(q^kq^H+q^{-k}q^{-H})^{-1}E=E(q^{k+1}q^H+q^{-k-1}q^{-H})^{-1}, \eqno (13) $$
$$
(q^kq^H+q^{-k}q^{-H})^{-1}F=F(q^{k-1}q^H+q^{-k+1}q^{-H})^{-1}. \eqno (14)
$$

Note that the algebra $U_q({\rm sl}_2)$ has finite dimensional irreducible
representations $T_l\equiv T_l^{(1)}$, $T_l^{(-1)}$, $T_l^{({\rm i})}$,
$T_l^{(-{\rm i})}$, $l=0,{\frac 12}, 1, {\frac 32}, \cdots$, acting on the
vector spaces ${\cal H}_l$ with bases $| m\rangle$, $m=-l,-l+1,\cdots ,l$.
These representations are given by the formulas
$$
T_l^{(1)} (q^H)|m\rangle =q^m |m\rangle ,\ \ \
T_l^{(1)} (E)|m\rangle =[l-m] |m+1\rangle , \eqno (15) $$
$$
T_l^{(1)} (F)|m\rangle =[l+m] |m-1\rangle , \eqno (16)
$$
where a number in square brackets means a $q$-number, defined by the formula
$$
[a]=\frac {q^a-q^{-a}}{q-q^{-1}} ,
$$
and by the formulas
$$
T_l^{(-1)} (q^H)|m\rangle =-q^m |m\rangle ,\ \ \
T_l^{(-1)} (E)=T_l^{(1)} (E),\ \ \ T_l^{(-1)} (F)=T_l^{(1)} (F), \eqno (17)
$$
$$
T_l^{({\rm i})} (q^H)|m\rangle ={\rm i}q^m |m\rangle ,\ \ \
T_l^{({\rm i})} (E)=T_l^{(1)} (E),\ \ \ T_l^{({\rm i})} (F)=-T_l^{(1)} (F),
\eqno (18)
$$
$$
T_l^{(-{\rm i})} (q^H)|m\rangle =-{\rm i}q^m |m\rangle ,\ \ \
T_l^{(-{\rm i})} (E)=T_l^{(1)} (E),\ \ \ T_l^{(-{\rm i})} (F)=-T_l^{(1)} (F).
\eqno (19)
$$
The representations $T_l^{(1)}$, $T_l^{(-1)}$,
$T_l^{({\rm i})}$, $T_l^{(-{\rm i})}$, $l=0,\frac 12 ,1,\frac 32 ,\cdots ,$
are pairwise non-equivalent, and any finite dimensional irreducible
representation of $U_q({\rm sl}_2)$ is equivalent to one of these
representations (see, for example, [11], Chapter 3).

Now we wish to extend these representations of $U_q({\rm sl}_2)$ to the
representations of ${\hat U}_q({\rm sl}_2)$ by using the
relation
$$
T((q^kq^H+q^{-k}q^{-H})^{-1}):=(q^kT(q^H)+q^{-k}T(q^{-H}))^{-1}.
$$
Clearly, only those irreducible representations $T$ of $U_q({\rm sl}_2)$
can be extended to ${\hat U}_q({\rm sl}_2)$ for which the operators
$q^kT(q^H)+q^{-k}T(q^{-H})$ are invertible. From formulas (15)--(19) it is
clear that these operators are always invertible for the representations
$T_l^{(1)}$, $T_l^{(-1)}$, $l=0,\frac 12 ,1,\frac 32 ,\cdots ,$ and for
the representations
$T_l^{({\rm i})}$, $T_l^{(-{\rm i})}$, $l=\frac 12 ,\frac 32 ,
\frac 52 ,\cdots $. For the representations
$T_l^{({\rm i})}$, $T_l^{(-{\rm i})}$, $l=0,1,2,
\cdots $, some of these operators are not invertible since they have
zero eigenvalue. Denoting the extended representations by the same symbols,
we can formulate the following statement:
\medskip

{\bf Proposition 1.} {\it The algebra ${\hat U}_q({\rm sl}_2)$
has the irreducible finite dimensional representations
$T_l^{(1)}$, $T_l^{(-1)}$, $l=0,\frac 12 ,1,\frac 32 ,\cdots ,$ and
$T_l^{({\rm i})}$, $T_l^{(-{\rm i})}$, $l=\frac 12 ,\frac 32 ,
\frac 52 ,\cdots $. Any irreducible finite dimensional representation
of ${\hat U}_q({\rm sl}_2)$ is equivalent to one of these representations.}
\bigskip

\noindent
{\sf III. THE ALGEBRA HOMOMORPHISM
$U_q({\rm so}_3)\to {\hat U}_q({\rm sl}_2)$}
\medskip

The aim of this section is to give (in an explicit form) the homomorphism
of the algebra $U_q({\rm so}_3)$ to ${\hat U}_q({\rm sl}_2)$. This
homomorphism is described by the following proposition:
\medskip

{\bf Proposition 2.} {\it There exists a unique algebra homomorphism
$\psi : U_q({\rm so}_3)\to {\hat U}_q({\rm sl}_2)$ such that
$$
\psi (I_1)=\frac {\rm i}{q-q^{-1}} (q^H-q^{-H}), \eqno (20) $$
$$
\psi (I_2)=(E-F) (q^H+q^{-H})^{-1}, \eqno (21) $$
$$
\psi (I_3)=({\rm i} q^{H-1/2}E+{\rm i}q^{-H-1/2}F) (q^H+q^{-H})^{-1},
\eqno (22)
$$
where $q^{H+a}:= q^Hq^a$ for $a\in {\bf C}$.}
\medskip

{\sl Proof.} In order to prove this proposition we have to show that
$$
q^{1/2}\psi (I_1)\psi (I_2)-q^{-1/2}\psi (I_2)\psi (I_1)=\psi (I_3), $$
$$
q^{1/2}\psi (I_2)\psi (I_3)-q^{-1/2}\psi (I_3)\psi (I_2)=\psi (I_1),
\eqno (23)$$
$$
q^{1/2}\psi (I_3)\psi (I_1)-q^{-1/2}\psi (I_1)\psi (I_3)=\psi (I_2).
$$
Let us prove the relation (23). (Other relations are proved similarly.)
Substituting the expressions (20)--(22) for
$\psi (I_i)$, $i=1,2,3$, into (23) we have (after multiplying both
sides of equality by $(q^H+q^{-H})$ on the right) the relation
$$
q(E-F)Eq^H(qq^H+q^{-1}q^{-H})^{-1}+q(E-F)Fq^{-H}(q^{-1}q^H+qq^{-H})^{-1}- $$
$$
-qE^2q^H(qq^H+q^{-1}q^{-H})^{-1}-q^{-1}FEq^{-H}(qq^H+q^{-1}q^{-H})^{-1}+ $$
$$
+q^{-1}EFq^H(q^{-1}q^H+qq^{-H})^{-1}+qF^2q^{-H}(q^{-1}q^H+qq^{-H})^{-1}
={\rm i}\frac {q^{2H}-q^{-2H}}{q-q^{-1}}.
$$
The formula (23) is true if and only if this relation is correct. We
multiply both its sides by $(qq^H+q^{-1}q^{-H})(q^{-1}q^H+qq^{-H})$
on the right and obtain the relation in the algebra $U_q({\rm sl}_2)$
(that is, without the expressions $(q^kq^H+q^{-k}q^{-H})^{-1}$). This
relation is easily verified by using the defining relations (9) and (10)
of the algebra $U_q({\rm sl}_2)$. Proposition is proved.
\bigskip

\noindent
{\sf IV. FINITE DIMENSIONAL REPRESENTATIONS OF $U_q({\rm so}_3)$:
$q$ IS NOT A ROOT OF UNITY}
\medskip

We assume in Sections IV--VII that $q$ is not a root of unity.

If $T$ is a representation of the algebra ${\hat U}_q({\rm sl}_2)$ on a
linear space ${\cal V}$, then the mapping $R: U_q({\rm so}_3)\to {\cal V}$
defined as the composition $R=T\circ \psi$, where $\psi$ is the
homomorphism from Proposition 2, is a representation of $U_q({\rm so}_3)$.
Let us consider the representations
$$
R_l^{(1)}=T_l^{(1)}\circ \psi ,\ \ \
R_l^{(-1)}=T_l^{(-1)}\circ \psi ,\ \ \
R_l^{({\rm i})}=T_l^{({\rm i})}\circ \psi ,\ \ \
R_l^{(-{\rm i})}=T_l^{(-{\rm i})}\circ \psi
$$
of $U_q({\rm so}_3)$, where $T_l^{(1)}$, $T_l^{(-1)}$,
$T_l^{({\rm i})}$, $T_l^{(-{\rm i})}$ are the irreducible representations
of ${\hat U}_q({\rm sl}_2)$ from Proposition 1.

Using formulas for the representations $T_l^{(\pm 1)}$ of $U_q({\rm sl}_2)$
and the expressions (20)--(22) for $\psi (I_j)$, $j=1,2,3$, we find that
$$
R_l^{(1)}(I_1)| m\rangle ={\rm i}[m]|m\rangle ,$$
$$
R_l^{(1)}(I_2)| m\rangle =\frac 1{q^m+q^{-m}} \{ [l-m] |m+1\rangle
-[l+m] |m-1\rangle \} , $$
$$
R_l^{(1)}(I_3)| m\rangle =\frac {{\rm i}q^{1/2}}{q^m+q^{-m}}
\{ q^m[l-m] |m+1\rangle
+q^{-m}[l+m] |m-1\rangle \}
$$
for the representation $R_l^{(1)}$ and
$$
R_l^{(-1)}(I_1)|m\rangle =-{\rm i}[m]|m\rangle ,\ \ \
R_l^{(-1)}(I_2)= -R_l^{(1)}(I_2),\ \ \  R_l^{(-1)}(I_3)=R_l^{(1)}(I_3).
$$
Denoting the vectors $|m\rangle$ by $|-m\rangle$ for the representations
$R_l^{(-1)}$ we easily find that the matrices of the representation
$R_l^{(-1)}$ in the basis $|-m\rangle$, $m=-l,-l+1,\cdots ,l$, coincide
with the corresponding matrices of the representation $R_l^{(1)}$. Thus,
the non-equivalent representations $T_l^{(1)}$ and $T_l^{(-1)}$ of the
algebra ${\hat U}_q({\rm sl}_2)$ lead to equivalent representations of
$U_q({\rm so}_3)$.

For the representations $R_l^{({\rm i})}$ and $R_l^{(-{\rm i})}$
we have
$$
R_l^{({\rm i})}(I_1) |m\rangle =-\frac {q^m+q^{-m}}{q-q^{-1}} |m\rangle ,$$
$$
R_l^{({\rm i})}(I_2) |m\rangle ={\rm i}\frac {[l-m]}{q^m-q^{-m}} |m+1\rangle
+{\rm i}\frac {[l+m]}{q^m-q^{-m}} |m-1\rangle ,$$
$$
R_l^{({\rm i})}(I_3) |m\rangle =-\frac {{\rm i}q^{m+1/2}[l-m]}{q^m-q^{-m}}
|m+1\rangle
-\frac {{\rm i}q^{-m+1/2}[l+m]}{q^m-q^{-m}} |m-1\rangle
$$
and
$$
R_l^{(-{\rm i})}(I_1) |m\rangle =\frac {q^m+q^{-m}}{q-q^{-1}} |m\rangle ,$$
$$
R_l^{(-{\rm i})}(I_2) |m\rangle =-{\rm i}\frac {[l-m]}{q^m-q^{-m}} |m+1\rangle
-{\rm i}\frac {[l+m]}{q^m-q^{-m}} |m-1\rangle ,$$
$$
R_l^{(-{\rm i})}(I_3) |m\rangle =-\frac {{\rm i}q^{m+1/2}[l-m]}{q^m-q^{-m}}
|m+1\rangle
-\frac {{\rm i}q^{-m+1/2}[l+m]}{q^m-q^{-m}} |m-1\rangle .
$$

{\bf Proposition 3.} {\it The representations $R_l^{(1)}$ of
$U_q({\rm so}_3)$ are irreducible. The representations
$R_l^{({\rm i})}$ and $R_l^{(-{\rm i})}$ are reducible.}
\medskip

{\sl Proof.} To prove the first part of the proposition we first note
that since $q$ is not a root of unity, the eigenvalues
${\rm i}[m]$, $m=-l,-l+1,\cdots ,l$, of the operator $R_l^{(1)}(I_1)$
are pairwise different.

Let $V$ be an invariant subspace of the space
${\cal H}_l$ of the representation $R_l^{(1)}$, and let ${\bf v}\equiv
\sum _{m_i} \alpha _i| m_i\rangle \in V$, where $| m_i\rangle $ are
eigenvectors of $R_l^{(1)}(I_1)$. Then $|m_i\rangle \in V$. We prove
this for the case when ${\bf v}=\alpha _1|m_1\rangle +\alpha _2|m_2\rangle$.
(The case of more number of summands is proved similarly.) We have
$R_l^{(1)}(I_1){\bf v}={\rm i}\alpha _1[m_1]|m_1\rangle +{\rm i}\alpha _2[m_2]
|m_2\rangle$. Since
$$
{\bf v}=\alpha _1|m_1\rangle +\alpha _2|m_2\rangle \in V,\ \ \
{\bf v}'\equiv {\rm i}\alpha _1[m_1]|m_1\rangle +{\rm i}\alpha _2[m_2]
|m_2\rangle \in V
$$
one derives that
$$
{\rm i}[m_1]{\bf v}-{\bf v}'={\rm i}\alpha _2([m_1]-[m_2])| m_2\rangle \in V.
$$
Since $[m_1]\ne [m_2]$, then  $|m_2\rangle \in V$ and hence
$|m_1\rangle \in V$.

In order to prove that $V={\cal H}_l$ we obtain from the above formulas for
$R_l^{(1)}(I_2)|m\rangle$ and $R_l^{(1)}(I_3)|m\rangle$ that
$$
\{ R_l^{(1)}(I_3)-{\rm i} q^{m+1/2}R_l^{(1)}(I_2)\} |m\rangle
={\rm i}q^{1/2}|m-1\rangle , $$
$$
\{ R_l^{(1)}(I_3)+{\rm i} q^{-m+1/2}R_l^{(1)}(I_2)\} |m\rangle
={\rm i}q^{1/2}|m+1\rangle .
$$
Since $V$ contains at least one basis vector $|m\rangle$, it follows
from these relations that $V$ contains the vectors $|m-1\rangle ,
|m-2\rangle ,\cdots ,|-l\rangle$ and the vectors
$|m+1\rangle ,|m+2\rangle ,\cdots ,|l\rangle$. This means that
$V={\cal H}_l$ and the representation $R_l^{(1)}$ is irreducible.

Let us show that the representations $R_l^{({\rm i})}$ are reducible.
The eigenvalues of the operator $R_l^{({\rm i})}(I_1)$ are
$$
-\frac {q^m+q^{-m}}{q-q^{-1}},\ \ \ m=-l,-l+1,\cdots ,l,
$$
that is, every spectral point has multiplicity 2. Namely, the pairs of vectors
$|m\rangle$ and $|-m\rangle$ are of the same eigenvalue.
Let $V_1$ be the subspace of the representation space
${\cal H}_l$ spanned by the vectors
$$\textstyle
|\frac 12 \rangle +{\rm i}\, |- \frac 12 \rangle ,\ \
|\frac 32 \rangle -{\rm i}\, |- \frac 32 \rangle ,\ \
|\frac 52 \rangle +{\rm i}\, |- \frac 52 \rangle ,\ \
|\frac 72 \rangle -{\rm i}\, |- \frac 72 \rangle , \ \cdots \ ,
\eqno (24)
$$
and let $V_2$ be the subspace spanned by the vectors
$$\textstyle
|\frac 12 \rangle -{\rm i}\, |- \frac 12 \rangle ,\ \
|\frac 32 \rangle +{\rm i}\, |- \frac 32 \rangle ,\ \
|\frac 52 \rangle -{\rm i}\, |- \frac 52 \rangle ,\ \
|\frac 72 \rangle +{\rm i}\, |- \frac 72 \rangle , \ \cdots \ .
\eqno (25)
$$
We denote the vectors (24) by
$$\textstyle
|\frac 12 \rangle ',\ \ \ |\frac 32 \rangle ',\ \ \
|\frac 52 \rangle ', \ \ \ |\frac 72 \rangle ' ,\ \cdots \eqno (26)
$$
and the vectors (25) by
$$\textstyle
|\frac 12 \rangle '',\ \ \ |\frac 32 \rangle '',\ \ \
|\frac 52 \rangle '', \ \ \ |\frac 72 \rangle '' ,\ \cdots . \eqno (27)
$$
Then
$$
R^{({\rm i})}_l (I_1)| m\rangle '=-\frac {q^m+q^{-m}}{q-q^{-1}}
| m\rangle ',\ \ \
R^{({\rm i})}_l (I_1)| m\rangle ''=-\frac {q^m+q^{-m}}{q-q^{-1}}
| m\rangle ''.
$$
We also have
$$
R^{({\rm i})}_l (I_2) {\textstyle  |\frac 12 \rangle '}=
{\rm i} \frac {[l-\frac 12]}{q^{1/2}-q^{-1/2}}{\textstyle |\frac 32 \rangle }
+{\rm i} \frac {[l+\frac 12]}{q^{1/2}-q^{-1/2}}
{\textstyle |-\frac 12 \rangle } +
\qquad\qquad\qquad\qquad\qquad   $$
$$\qquad\qquad\qquad\qquad\qquad
+\frac {[l+\frac 12]}{q^{1/2}-q^{-1/2}}{\textstyle |\frac 12 \rangle }
+ \frac {[l-\frac 12]}{q^{1/2}-q^{-1/2}}{\textstyle |-\frac 32 \rangle }  $$
$$
=\frac {[l+\frac 12]}{q^{1/2}-q^{-1/2}}{\textstyle |\frac 12 \rangle '}
+{\rm i} \frac {[l-\frac 12]}{q^{1/2}-q^{-1/2}}
{\textstyle |\frac 32 \rangle '} .
$$
We derive similarly that
$$
R^{({\rm i})}_l (I_2){\textstyle  |\frac 12 \rangle ''}=
-\frac {[l+\frac 12]}{q^{1/2}-q^{-1/2}}{\textstyle |\frac 12 \rangle ''}
+{\rm i} \frac {[l-\frac 12]}{q^{1/2}-q^{-1/2}}
{\textstyle |\frac 32 \rangle ''}
$$
and that
$$
R^{({\rm i})}_l (I_2) |m \rangle '=
{\rm i}\frac {[l-m]}{q^{m}-q^{-m}}| m+1 \rangle '
+{\rm i} \frac {[l+m]}{q^{m}-q^{-m}}|m-1 \rangle ', \ \ \
m>{\textstyle \frac 12} ,
$$
$$
R^{({\rm i})}_l (I_2) |m \rangle ''=
{\rm i}\frac {[l-m]}{q^{m}-q^{-m}}| m+1 \rangle ''
+{\rm i} \frac {[l+m]}{q^{m}-q^{-m}}|m-1 \rangle '',\ \ \
m>{\textstyle \frac 12} .
$$
Thus, the subspaces $V_1$ and $V_2$ are invariant with respect to the
operators $R^{({\rm i})}_l (I_1)$ and $R^{({\rm i})}_l (I_2)$. This means
that they are invariant with respect to the representation
$R^{({\rm i})}_l$.

It is proved similarly that the subspace $V_1$ of the space ${\cal H}_l$ of
the representation $R^{(-{\rm i})}_l$ spanned by the vectors (24) and the
subspace $V_2$ of ${\cal H}_l$ spanned by the vectors (25) are invariant
with respect
to the operators $R^{(-{\rm i})}_l (I_1)$ and $R^{(-{\rm i})}_l (I_2)$.
That is, the representation $R^{(-{\rm i})}_l$ is also reducible.
Proposition is proved.
\medskip

Let $R^{({\rm i},+)}_n$ and $R^{({\rm i},-)}_n$, $n=l+
\frac 12 ={\rm dim}\, V_1
={\rm dim}\, V_2$, be the representations of $U_q({\rm so}_3)$ which
are restrictions of $R^{({\rm i})}_l$ to the
subspaces $V_1$ and $V_2$, respectively. Denoting the vectors (26) of
the subspace $V_1$ by
$$\textstyle
|1\rangle ,\ \ \  |2\rangle ,\ \ \ |3\rangle ,\ \ \ |4\rangle ,\ \cdots ,\
|n\rangle \equiv | l+\frac 12 \rangle , \eqno (28)
$$
respectively, we have
$$
R^{({\rm i},+)}_n (I_1)| k\rangle =-\frac {q^{k-1/2}+q^{-k+1/2}}{q-q^{-1}}
| k\rangle , $$
$$
R^{({\rm i},+)}_n (I_2)| 1\rangle =\frac {[n]}{q^{1/2}-q^{-1/2}}
| 1\rangle
+{\rm i} \frac {[n-1]}{q^{1/2}-q^{-1/2}}
| 2\rangle , $$
$$
R^{({\rm i},+)}_n (I_2)| k\rangle ={\rm i}\frac {[n-k]}{q^{k-1/2}-
q^{-k+1/2}} | k+1\rangle +
{\rm i} \frac {[n+k-1]}{q^{k-1/2}-q^{-k+1/2}}
| k-1\rangle , \ \ \ k\ne 1.
$$
For the operator $R^{({\rm i},+)}_n (I_3)$ we have
$$
R^{({\rm i},+)}_n (I_3)| 1\rangle =-\frac {[n]}{q^{1/2}-q^{-1/2}}
| 1\rangle
-{\rm i} \frac {q[n-1]}{q^{1/2}-q^{-1/2}} | 2\rangle , $$
$$
R^{({\rm i},+)}_n (I_3)| k\rangle =-{\rm i}\frac {q^k[n-k]}
{q^{k-1/2}-q^{-k+1/2}} | k+1\rangle
-{\rm i} \frac {q^{-k+1}[n+k-1]}{q^{k-1/2}-q^{-k+1/2}} | k-1\rangle ,\ \
k\ne 1.
$$

Denoting the vectors (27) of the subspace $V_2$ by the symbols (28),
respectively, we obtain
$$
R^{({\rm i},-)}_n (I_1)| k\rangle =-\frac {q^{k-1/2}+q^{-k+1/2}}{q-q^{-1}}
| k\rangle , $$
$$
R^{({\rm i},-)}_n (I_2)| 1\rangle =-\frac {[n]}{q^{1/2}-q^{-1/2}}
| 1\rangle
+{\rm i} \frac {[n-1]}{q^{1/2}-q^{-1/2}}
| 2\rangle , $$
$$
R^{({\rm i},-)}_n (I_2)| k\rangle
=R^{({\rm i},+)}_n (I_2)| k\rangle ,  \ \ \ k\ne 1.
$$
For the operator $R^{({\rm i},-)}_l (I_3)$ we find that
$$
R^{({\rm i},-)}_n (I_3)| 1\rangle =\frac {[n]}{q^{1/2}-q^{-1/2}}
| 1\rangle
-{\rm i} \frac {q[n-1]}{q^{1/2}-q^{-1/2}} | 2\rangle , $$
$$
R^{({\rm i},-)}_n (I_3)| k\rangle =
R^{({\rm i},+)}_n (I_3)| k\rangle , \ \ \ k\ne 1.
$$

Let now  $R^{(-{\rm i},+)}_n$ and $R^{(-{\rm i},-)}_n$,
$n=l+\frac 12 $, be the representations of $U_q({\rm so}_3)$ which
are restrictions of the representation $R^{(-{\rm i})}_l$ to the
subspaces $V_1$ and $V_2$, respectively. Introducing the vectors similar
to the vectors (28), for the representation $R^{(-{\rm i},+)}_n$
we have
$$
R^{(-{\rm i},+)}_n (I_1)| k\rangle =\frac {q^{k-1/2}+q^{-k+1/2}}{q-q^{-1}}
| k\rangle , \ \ \
R^{(-{\rm i},+)}_n (I_2)= - R^{({\rm i},+)}_n (I_2), $$
$$
R^{(-{\rm i},+)}_n (I_3)| 1\rangle =\frac {[n]}{q^{1/2}-q^{-1/2}}
| 1\rangle
+{\rm i} \frac {q[n-1]}{q^{1/2}-q^{-1/2}}
| 2\rangle , $$
$$
R^{(-{\rm i},+)}_n (I_3)| k\rangle ={\rm i}\frac {q^k[n-k]}{q^{k-1/2}-
q^{-k+1/2}} | k+1\rangle +
{\rm i} \frac {q^{-k+1}[n+k-1]}{q^{k-1/2}-q^{-k+1/2}}
| k-1\rangle , \ \ \ k\ne 1.
$$
For the representation $R^{(-{\rm i},-)}_l (I_3)$ we obtain
$$
R^{(-{\rm i},-)}_n (I_1)| k\rangle =\frac {q^{k-1/2}+q^{-k+1/2}}{q-q^{-1}}
| k\rangle , \ \ \
R^{(-{\rm i},-)}_n (I_2)= - R^{({\rm i},-)}_n (I_2), $$
$$
R^{(-{\rm i},-)}_n (I_3)| 1\rangle =-\frac {[n]}{q^{1/2}-q^{-1/2}}
| 1\rangle
+{\rm i} \frac {q[n-1]}{q^{1/2}-q^{-1/2}}
| 2\rangle , $$
$$
R^{(-{\rm i},-)}_n (I_3)| k\rangle =
R^{(-{\rm i},+)}_n (I_3)| k\rangle .
$$

Thus, we constructed the representations
$R^{({\rm i},+)}_n$,
$R^{({\rm i},-)}_n$, $R^{(-{\rm i},+)}_n$ and $R^{(-{\rm i},-)}_n$
of the algebra $U_q({\rm so}_3)$. The following theorem characterizes
them.
\medskip

{\bf Theorem 1.} {\it The representations $R^{({\rm i},+)}_n$,
$R^{({\rm i},-)}_n$, $R^{(-{\rm i},+)}_n$ and $R^{(-{\rm i},-)}_n$
are irreducible and pairwise nonequivalent. For any $l$ the representation
$R^{(1)}_l$ is not equivalent to any of these representations.}
\medskip

{\sl Proof.} The irreducibility is proved exactly in the same way
as in Proposition 3. Equivalence relations may exist only
for irreducible representations of the same dimension. That is, we
have to show that under fixed $n$ no pair of the representations
$R^{({\rm i},+)}_n$,
$R^{({\rm i},-)}_n$, $R^{(-{\rm i},+)}_n$ and $R^{(-{\rm i},-)}_n$
is equivalent. It follows from the above formulas that
the operators $R^{({\rm i},+)}_n (I_1)$ and $R^{({\rm i},-)}_n (I_1)$,
as well as the operators
$R^{(-{\rm i},+)}_n (I_1)$ and $R^{(-{\rm i},-)}_n (I_1)$,
have the same set of eigenvalues. Moreover, the spectrum of the first
pair of operators differs from that of the second pair.
Hence, no of representations $R^{({\rm i},+)}_n $ and $R^{({\rm i},-)}_n$
is equivalent to $R^{(-{\rm i},+)}_n $ or $R^{(-{\rm i},-)}_n$.
The representations $R^{({\rm i},+)}_n $ and $R^{({\rm i},-)}_n$
are not equivalent since the operators
$R^{({\rm i},+)}_n (I_2) $ and $R^{({\rm i},-)}_n (I_2)$
have different traces (for equivalent representations these operators
must have the same trace). For the same reason, the representations
$R^{(-{\rm i},+)}_n $ and $R^{(-{\rm i},-)}_n$ are not equivalent.
The last assertion of the theorem follows from the fact that the spectrum of
the operator $R^{(1)}_l(I_1)$ differs from the spectra of the operators
$R^{({\rm i},+)}_n(I_1)$, $R^{({\rm i},-)}_n(I_1)$,
$R^{(-{\rm i},+)}_n(I_1)$ and $R^{(-{\rm i},-)}_n(I_1)$.
Theorem is proved.
\medskip

Clearly, the reducible representations $R^{({\rm i})}_n$ and
$R^{(-{\rm i})}_n$ decomposes into irreducible components as
$$
R^{({\rm i})}_n= R^{({\rm i},+)}_n\oplus  R^{({\rm i},-)}_n,\ \ \
R^{(-{\rm i})}_n= R^{(-{\rm i},+)}_n\oplus  R^{(-{\rm i},-)}_n.
\eqno (29)
$$

It can be proved that {\it every irreducible finite dimensional
representation of $U_q({\rm so}_3)$ is equivalent to one of the
representations $R^{(1)}_l$,
$R^{({\rm i},+)}_n$, $R^{({\rm i},-)}_n$,
$R^{(-{\rm i},+)}_n$, $R^{(-{\rm i},-)}_n$}. That is, these representations
exhaust, up to equivalence, all irreducible finite dimensional
representations of $U_q({\rm so}_3)$. A proof of this statement will be given
in a separate paper.
\bigskip

\noindent
{\sf V. TENSOR PRODUCTS OF REPRESENTATIONS OF $U_q({\rm so}_3)$}
\medskip

As mentioned above, no Hopf algebra structure is known for the algebra
$U_q({\rm so}_3)$.
Therefore, we cannot construct tensor product of finite
dimensional representations of $U_q({\rm so}_3)$ by using a comultiplication
as we do in the case of the quantum algebra $U_q({\rm sl}_2)$.
However, we may construct some tensor product representations by using
the algebra homomorphism of Proposition 2.

First we determine which tensor products of irreducible representations of
$U_q({\rm sl}_2)$ can be extended to representations of the algebra
${\hat U}_q({\rm sl}_2)$. Verifying for which tensor products $T=T'\otimes
T''$ of irreducible representations of $U_q({\rm sl}_2)$ the operators
$$
q^kT(q^H)+q^{-k}T(q^{-H}),\ \ \ \ k\in {\bf Z},
$$
are invertible, we conclude that only the tensor products
$$\textstyle
T^{(\pm 1)}_l\otimes T^{(\pm 1)}_{l'},\ \
T^{(\pm 1)}_l\otimes T^{(\mp 1)}_{l'},\ \ \ l,l'=0,\frac 12 ,1, \frac 32 ,
\cdots , $$
$$\textstyle
T^{(\pm 1)}_l\otimes T^{(\pm {\rm i})}_{l'},\ \
T^{(\pm 1)}_l\otimes T^{(\mp {\rm i})}_{l'},\ \ \
l=0,1,2,\cdots ,\ \ l'=\frac 12 ,\frac 32 ,\frac 52 ,\cdots ,
$$
$$\textstyle
T^{(\pm {\rm i})}_l\otimes T^{(\pm 1)}_{l'},\ \
T^{(\pm {\rm i})}_l\otimes T^{(\mp 1)}_{l'},\ \ \
l=\frac 12 ,\frac 32 ,\frac 52 ,\cdots ,\ \  l'=0,1,2,\cdots ,
$$
$$\textstyle
T^{(\pm {\rm i})}_l\otimes T^{(\pm {\rm i})}_{l'},\ \
T^{(\pm {\rm i})}_l\otimes T^{(\mp {\rm i})}_{l'},\ \
l, l'=\frac 12 ,\frac 32 ,\frac 52 ,\cdots ,
$$
can be extended to the algebra ${\hat U}_q({\rm sl}_2)$. Taking into account
the decompositions of tensor products of irreducible representations of
$U_q({\rm sl}_2)$ (see, for example, the end of Subsection 3.2.1 and
Proposition 3.22 in [11]) we find that
$$
T^{(\omega )}_l\otimes T^{(\omega ')}_{l'} \simeq
T^{(\omega \omega ')}_{l+l'} \oplus
T^{(\omega \omega ')}_{l+l'-1} \oplus \cdots \oplus
T^{(\omega \omega ')}_{|l+l'|} , \eqno (30) $$
$$
T^{(\omega )}_l\otimes T^{(\pm {\rm i})}_{l'} \simeq
T^{(\pm \omega {\rm i})}_{l+l'} \oplus
T^{(\pm \omega {\rm i})}_{l+l'-1} \oplus \cdots \oplus
T^{(\pm \omega {\rm i})}_{|l+l'|} , \eqno (31) $$
$$
T^{(\pm {\rm i})}_{l} \otimes T^{(\omega )}_{l'} \simeq
T^{(\pm \omega {\rm i})}_{l+l'} \oplus
T^{(\pm \omega {\rm i})}_{l+l'-1} \oplus \cdots \oplus
T^{(\pm \omega {\rm i})}_{|l+l'|} , \eqno (32) $$
$$
T^{(\omega {\rm i})}_{l} \otimes T^{(\omega '{\rm i})}_{l'} \simeq
T^{(- \omega \omega ')}_{l+l'} \oplus
T^{(- \omega \omega ')}_{l+l'-1} \oplus \cdots \oplus
T^{(- \omega \omega ')}_{|l+l'|} , \eqno (33)
$$
where $\omega ,\omega ' =\pm 1$.

Now we define tensor products of representations of
$U_q({\rm so}_3)$ corresponding to the above tensor product representations
of ${\hat U}_q({\rm sl}_2)$ as
$$
R\otimes R'=(T\otimes T')\circ \psi ,
$$
where $R=T\circ \psi$ and $R'=T'\circ \psi$. Taking into account the
definitions of tensor products of representations of $U_q({\rm sl}_2)$
by means of the comultiplication and the definition of the mapping
$\psi$ we have
$$
(R\otimes R')(I_1)=(T\otimes T')\circ \psi (I_1)$$
$$
=\frac {\rm i}{q-q^{-1}}
\left( T(q^H)\otimes T'(q^H) - T(q^{-H})\otimes T'(q^{-H})\right) .
$$
Similarly,
$$
(R\otimes R')(I_2)=\left( T(E)\otimes T'(q^H)+T(q^{-H})\otimes T'(E) -
T(F)\otimes T'(q^H)- \right. $$
$$\left.
-(T(q^{-H})\otimes T'(F)\right)
\left( T(q^H)\otimes T'(q^H) + T(q^{-H})\otimes T'(q^{-H})\right) ^{-1} .
$$

Composing both sides of the relations (30)--(33) with the mapping $\psi$
of Proposition 2, we find the decomposition into representations of
$U_q({\rm so}_3)$ for the tensor products
$$
R^{(1)}_l\otimes R^{(1)}_{l'},\ \ \
R^{(1)}_l\otimes R^{(\pm {\rm i})}_{l'},\ \ \
R^{(\pm {\rm i})}_{l'}\otimes R^{(1)}_l,\ \ \
R^{(\pm {\rm i})}_l\otimes R^{(\pm {\rm i})}_{l'},\ \ \
R^{(\pm {\rm i})}_l\otimes R^{(\mp {\rm i})}_{l'},
$$
where the second and the third tensor products are defined only for
$l=0,1,2,\cdots $. (Note that the representations $R^{(\pm {\rm i})}_l$
are defined only for $l=\frac 12 ,\frac 32 ,\frac 52 ,\cdots$.) We have
$$
R^{(1)}_l\otimes R^{(1)}_{l'} \simeq
R^{(1)}_{l+l'} \oplus
R^{(1)}_{l+l'-1} \oplus \cdots \oplus
R^{(1)}_{|l-l'|} ,  $$
$$
R^{(1)}_l\otimes R^{(\pm {\rm i})}_{l'} \simeq
R^{(\pm {\rm i})}_{l+l'} \oplus
R^{(\pm {\rm i})}_{l+l'-1} \oplus \cdots \oplus
R^{(\pm {\rm i})}_{|l-l'|} ,  $$
$$
R^{(\pm {\rm i})}_l\otimes R^{(1)}_{l'} \simeq
R^{(\pm {\rm i})}_{l+l'} \oplus
R^{(\pm {\rm i})}_{l+l'-1} \oplus \cdots \oplus
R^{(\pm {\rm i})}_{|l-l'|} ,  $$
$$
R^{(\omega {\rm i})}_l\otimes R^{(\omega '{\rm i})}_{l'} \simeq
R^{(1)}_{l+l'} \oplus
R^{(1)}_{l+l'-1} \oplus \cdots \oplus
R^{(1)}_{|l-l'|} .
$$
In these formulas the representations $R^{(\pm {\rm i})}_l$ are reducible.
Unfortunately, our definition of tensor products of representations of
$U_q({\rm so}_3)$ does not allow to determine the tensor products
containing the irreducible
representations $R^{(\pm {\rm i},\pm )}_n$ and
$R^{(\pm {\rm i},\mp )}_n$.
\bigskip

\noindent
{\sf VI. INFINITE DIMENSIONAL REPRESENTATIONS OF $U_q({\rm so}_3)$
OBTAINED FROM REPRESENTATIONS OF $U_q({\rm sl}_2)$}
\medskip

By using the homomorphism $\psi : U_q({\rm so}_3)\to {\hat U}_q({\rm sl}_2)$
from Proposition 2 and infinite dimensional irreducible representations
of the algebra ${\hat U}_q({\rm sl}_2)$ we can construct infinite
dimensional irreducible representations of the algebra $U_q({\rm so}_3)$.

Let us first describe irreducible infinite dimensional representations
of the algebra $U_q({\rm sl}_2)$. Note that by an infinite dimensional
representation $T$ of $U_q({\rm sl}_2)$ we mean a homomorphism of
$U_q({\rm sl}_2)$ into the algebra of linear operators (bounded or
unbounded) on a Hilbert space, defined on an everywhere dense invariant
subspace $D$, such that the operator $T(q^H)$ can be diagonalized, has a
discrete spectrum and its eigenvectors belong to $D$. Infinite dimensional
representations $T$ of $U_q({\rm so}_3)$ are described in the same way
replacing the operator $T(q^H)$ by $T(I_1)$.

Two representations $T$ and $T'$ of $U_q({\rm sl}_2)$ on spaces
${\cal H}$ and ${\cal H}'$, respectively, are called (algebraically)
equivalent if there exist everywhere dence invariant subspaces
$V\subset {\cal H}$ and $V'\subset {\cal H}'$ and a one-to-one linear
operator $A: V\to V'$ such that $AT(a){\bf v}=T'(a)A{\bf v}$ for all
$a\in U_q({\rm sl}_2)$ and ${\bf v}\in V$. Equivalence of infinite
dimensional representations of $U_q({\rm so}_3)$ is defined in the same way.

Let $\epsilon$ be a fixed complex number such that $0\le {\rm Re}\,
\epsilon < 1$, and let ${\cal H}_\epsilon$ be a complex Hilbert space
with the orthonormal basis
$$
|m\rangle ,\ \ \ \ \ m=n+\epsilon ,\ \ n=0,\pm 1,\pm 2,\cdots . \eqno (34)
$$
For every complex number $a$ we construct the representation
$T_{a\epsilon}$ on the Hilbert space ${\cal H}_\epsilon$ defined by
$$
T_{a\epsilon}(q^H)|m\rangle =q^m |m\rangle , \
T_{a\epsilon}(E)|m\rangle =[a-m]|m+1\rangle , \
T_{a\epsilon}(F)|m\rangle =[a+m]|m-1\rangle ,
$$
where $[a\pm m]$ is the $q$-number (see, for example, [12]).
The equivalence relations in the set of the representations $T_{a\epsilon}$
can be extracted from the paper [12].

Note that the representation $T_{a\epsilon}$ is irreducible
if and only if $a\ne \pm \epsilon \, ({\rm mod} \,{\bf Z})$.

All the representations $T_{a\epsilon}$ can be extended to representations
of the algebra ${\hat U}_q({\rm sl}_2)$ except for the case when
$\epsilon =\pm {\rm i}\pi /2\tau$, where $q=e^\tau$. (We suppose below that
$\epsilon \ne \pm {\rm i}\pi /2\tau$.) We denote these extended
representations by the same symbols $T_{a\epsilon}$.

The formula $R_{a\epsilon}=T_{a\epsilon}\circ \psi$ associates with every
irreducible representation $T_{a\epsilon}$,
$\epsilon \ne \pm {\rm i}\pi /2\tau $, of ${\hat U}_q({\rm sl}_2)$
a representation of the algebra $U_q({\rm so}_3)$.

Let $\epsilon \ne \pm {\rm i}\pi /2\tau$ and
$\epsilon \ne \pm {\rm i}\pi /2\tau +\frac 12$. Then for the representations
$R_{a\epsilon}$ of $U_q({\rm so}_3)$ we have
$$
R_{a\epsilon}(I_1)|m\rangle ={\rm i}[m] |m\rangle , \eqno (35a)$$
$$
R_{a\epsilon}(I_2)|m\rangle =\frac 1{q^m+q^{-m}} \left\{ [a-m]|m+1\rangle
-[a+m]|m-1\rangle \right\} , \eqno (35b)$$
$$
R_{a\epsilon}(I_3)|m\rangle =\frac {{\rm i}q^{1/2}}{q^m+q^{-m}}
\left\{q^m [a-m]|m+1\rangle
+q^{-m}[a+m]|m-1\rangle \right\} .  \eqno (35c)
$$
If $\epsilon = {\rm i}\pi /2\tau +\frac 12$, then denoting the basis
elements $|m\rangle$, $m=n+\epsilon$, $n\in {\bf Z}$, by $|n+\frac 12
\rangle$, $n\in {\bf Z}$, respectively, we obtain
$$
R_{a\epsilon}(I_1)|k\rangle =-\frac {q^k+q^{-k}}{q-q^{-1}} |k\rangle ,$$
$$
R_{a\epsilon}(I_2)|k\rangle ={\rm i}\frac {[a'-k]}{q^k-q^{-k}}|k+1\rangle
+{\rm i}\frac {[a'+k]}{q^k-q^{-k}} |k-1\rangle  ,$$
$$
R_{a\epsilon}(I_3)|k\rangle =-\frac {{\rm i}q^{k+1/2}[a'-k]}{q^k-q^{-k}}
|k+1\rangle -\frac {{\rm i}q^{-k+1/2}[a'+k]}{q^k-q^{-k}}
|k-1\rangle ,
$$
where $a'=a+{\rm i}\pi /2\tau$ and $k=n+\frac 12$. If
$\epsilon =- {\rm i}\pi /2\tau +\frac 12$, then using the same notations
for basis elements we obtain
$$
R'_{a\epsilon}(I_1)|k\rangle =\frac {q^k+q^{-k}}{q-q^{-1}} |k\rangle ,$$
$$
R'_{a\epsilon}(I_2)|k\rangle =-{\rm i}\frac {[a'-k]}{q^k-q^{-k}}|k+1\rangle
-{\rm i}\frac {[a'+k]}{q^k-q^{-k}} |k-1\rangle  ,$$
$$
R'_{a\epsilon}(I_3)|k\rangle =-\frac {{\rm i}q^{k+1/2}[a'-k]}{q^k-q^{-k}}
|k+1\rangle -\frac {{\rm i}q^{-k+1/2}[a'+k]}{q^k-q^{-k}}
|k-1\rangle
$$
(to distinquish these representations from the previous ones we supplied
$R_{a\epsilon}$ by prime).
\medskip

{\bf Proposition 4.} {\it The representations $R_{a\epsilon}$ of
$U_q({\rm so}_3)$ are irreducible for irreducible representations
$T_{a\epsilon}$, $\epsilon \ne \pm {\rm i}\pi /2\tau +\frac 12$,
of ${\hat U}_q({\rm sl}_2)$. The representations
$R_{a\epsilon}$, $\epsilon = {\rm i}\pi /2\tau +\frac 12$, and
$R'_{a\epsilon}$, $\epsilon =-{\rm i}\pi /2\tau +\frac 12$,
are reducible.}
\medskip

{\sl Proof} is given in the same way as in the case of Proposition 3.
\medskip

As in the case of finite dimensional representations in Section IV,
decomposing the representations
$R_{a\epsilon}$, $\epsilon =\pm {\rm i}\pi /2\tau +\frac 12$, we obtain
irreducible infinite dimensional representations
of $U_q({\rm so}_3)$ which will be denoted by $R^{({\rm i},\pm )}_{a'}$ and
$R^{(-{\rm i},\pm )}_{a'}$, $a'=a+{\rm i}\pi /2\tau$. In the basis
$$
|n\rangle ,\ \ \ n=1,2,3,\cdots ,
$$
they are given by the formulas
$$
R^{({\rm i},\pm )}_{a'} (I_1)| k\rangle =
-\frac {q^{k-1/2}+q^{-k+1/2}}{q-q^{-1}} | k\rangle , $$
$$
R^{({\rm i},\pm )}_{a'} (I_2)| 1\rangle =\pm \frac {[a']}{q^{1/2}-q^{-1/2}}
| 1\rangle
+{\rm i} \frac {[a'-1]}{q^{1/2}-q^{-1/2}}
| 2\rangle , $$
$$
R^{({\rm i},\pm )}_{a'} (I_2)| k\rangle ={\rm i}\frac {[a'-k]}{q^{k-1/2}-
q^{-k+1/2}} | k+1\rangle +
{\rm i} \frac {[a'+k-1]}{q^{k-1/2}-q^{-k+1/2}}
| k-1\rangle , \ \ \ k\ne 1.
$$
$$
R^{({\rm i},\pm )}_{a'} (I_3)| 1\rangle =\mp \frac {[a']}{q^{1/2}-q^{-1/2}}
| 1\rangle
-{\rm i} \frac {q[a'-1]}{q^{1/2}-q^{-1/2}} | 2\rangle , $$
$$
R^{({\rm i},\pm )}_{a'} (I_3)| k\rangle =-{\rm i}\frac {q^k[a'-k]}
{q^{k-1/2}-q^{-k+1/2}} | k+1\rangle
-{\rm i} \frac {q^{-k+1}[a'+k-1]}{q^{k-1/2}-q^{-k+1/2}} | k-1\rangle ,\ \
k\ne 1.
$$
and by the formulas
$$
R^{(-{\rm i},\pm )}_{a'} (I_1)| k\rangle =
\frac {q^{k-1/2}+q^{-k+1/2}}{q-q^{-1}} | k\rangle , \ \ \
R^{(-{\rm i},\pm )}_{a'} (I_2)= - R^{({\rm i},\pm )}_{a'} (I_2), $$
$$
R^{(-{\rm i},\pm )}_{a'} (I_3)| 1\rangle =\pm \frac {[a']}{q^{1/2}-q^{-1/2}}
| 1\rangle
+{\rm i} \frac {q[a'-1]}{q^{1/2}-q^{-1/2}} | 2\rangle , $$
$$
R^{(-{\rm i},\pm )}_{a'} (I_3)| k\rangle ={\rm i}\frac {q^k[a'-k]}{q^{k-1/2}-
q^{-k+1/2}} | k+1\rangle +
{\rm i} \frac {q^{-k+1}[a'+k-1]}{q^{k-1/2}-q^{-k+1/2}}
| k-1\rangle , \ \ \ k\ne 1.
$$

{\bf Theorem 2.} {\it The representations $R^{({\rm i},\pm )}_{a'}$
$R^{(-{\rm i},\pm )}_{a'}$ are irreducible and pairwise nonequivalent.
For any $a$ the irreducible representation $R_{a\epsilon}$ is not equivalent
to some of these representations.}
\medskip

{\sl Proof} is given in the same way as in the finite dimensional case
(see the proof of Theorem 1).
\medskip

The algebra $U_q({\rm sl}_2)$ has also irreducible infinite dimensional
representations with highest weights or with lowest weights. They are
classified in the paper [12]. All of these representations $T$ can be extended
to the algebra ${\hat U}_q({\rm sl}_2)$. Using the composition
$R=T\circ \psi$ we obtain the corresponding representations $R$ of
$U_q({\rm so}_3)$. As above, it can be easily proved that to
nonequivalent representations $T$ of ${\hat U}_q({\rm sl}_2)$
with highest or lowest weight
there correspond nonequivalent irreducible representations of
$U_q({\rm so}_3)$. We give a list of these representations.

Let $l=\frac 12 ,1, \frac 32 ,2,\cdots$. We denote by $R^+_l$ the
representation of $U_q({\rm so}_3)$ acting on the Hilbert space ${\cal H}_l$
with the orthonormal basis $|m\rangle$, $m=l,l+1,l+2,\cdots$, and given by
formulas (35) with $a=-l$. By $R^-_l$ we denote the
representation of $U_q({\rm so}_3)$ acting on the Hilbert space
${\hat{\cal H}}_l$ with the orthonormal basis $|m\rangle$,
$m=-l,-l-1,-l-2,\cdots$, and given by
formulas (35) with $a=l$.

Now let $a\ne 0$ (mod ${\bf Z}$) and $a\ne \frac 12$ (mod ${\bf Z}$). We
denote by ${\cal H}_a$ the Hilbert space
with the orthonormal basis $|m\rangle$, $m=-a,-a+1,-a+2,\cdots$.
On this space the representation $R^+_a$ acts which is given by formulas
(35). On the Hilbert space ${\hat{\cal H}}_a$
with the orthonormal basis $|m\rangle$, $m=a,a-1,a-2,\cdots$,
the representation $R^-_a$ acts which is given by formulas (35).
\medskip

{\bf Proposition 5.} {\it The above representations $R^\pm _l$ and $R^\pm _a$
are irreducible and pairwise nonequivalent.}
\medskip

Proof of this proposition is contained in [13].
\bigskip

\noindent
{\sf VII. OTHER INFINITE DIMENSIONAL REPRESENTATIONS OF
$U_q({\rm so}_3)$}
\medskip

The algebra $U_q({\rm so}_3)$ has also irreducible infinite dimensional
representations which cannot be obtained from representations of
${\hat U}_q({\rm sl}_2)$. We describe these representations in this section.

Let ${\cal H}$ be the infinite dimensional vector space with the basis
$|m\rangle$, $m=0,\pm1,\pm2$, $\cdots$, and let $\lambda =q^\tau$ be a
nonzero complex
number such that $0\le {\rm Re}\, \tau < 1$. Then a direct calculation
shows that the operators $Q^+_\lambda (I_1)$ and $Q^+_\lambda (I_2)$ given by
the formulas
$$
Q^+_\lambda (I_1)|m\rangle =\frac {\lambda q^m +\lambda ^{-1}q^{-m}}
{q-q^{-1}}\, |m\rangle , $$
$$
Q^+_\lambda (I_2)|m\rangle =\frac 1{q-q^{-1}}\, |m+1\rangle  +
\frac 1{q-q^{-1}}\, |m-1\rangle
$$
satisfy the relations (7) and (8) and hence determine a representation of
$U_q({\rm so}_3)$ which will be denoted by $Q^+_\lambda$. Similarly, the
operators $Q^-_\lambda (I_1)$ and $Q^-_\lambda (I_2)$ given on
the space ${\cal H}$ by
$$
Q^-_\lambda (I_1)|m\rangle =-\frac {\lambda q^m +\lambda ^{-1}q^{-m}}
{q-q^{-1}}\, |m\rangle ,\ \ \
Q^-_\lambda (I_2):= Q^+_\lambda (I_2)
$$
determine a representation of $U_q({\rm so}_3)$ which is denoted by
$Q^-_\lambda$. The operators $Q_\lambda ^{\pm}(I_3)$ can be calculated
by means of formula (4).
\medskip

{\bf Proposition 6.} {\it If $\lambda \ne 1$ and $\lambda \ne q^{1/2}$, then
the representations $Q^+_\lambda$ and $Q^-_\lambda$
are irreducible. The representations $Q^\pm _1$ and
$Q^\pm _{\sqrt q}$ are reducible.}
\medskip

{\sl Proof.} The first part is proved in the same way as that of
Proposition 3. Let us prove the second part. The representations $Q^\pm _1$
and $Q^\pm _{\sqrt q}$ are the only representations in the set
$\{ Q^\pm _\lambda\}$ for which
the operator $Q^\pm _\lambda (I_1)$ has not a simple spectrum. The operators
$Q^\pm _1(I_1)$ has the spectrum
$$
\cdots ,\ q^{-2}+q^2,\ q^{-1}+q,\ 2,\ q+q^{-1},\ q^2+q^{-2},\ \cdots .
$$
Thus, only the spectral point 2 has multiplicity 1. All other points have
multiplicity 2.
Let $V_1$ and $V_2$ be the vector subspaces of ${\cal H}$ with the bases
$$
|0\rangle ,\ \ |m\rangle '=|m\rangle -|-m\rangle ,\ \ m=1,2,\cdots ,
$$
and
$$
|m\rangle ''=|m\rangle +|-m\rangle ,\ \ m=1,2,\cdots ,
$$
respectively. These basis vectors are eigenvectors of the operator
$Q^\pm _1(I_1)$:
$$
Q^\pm _1(I_1)| m\rangle '=\pm \frac {q^m+q^{-m}}{q-q^{-1}} \,
|m\rangle ',
\ \ \ Q^\pm _1(I_1)| m\rangle ''=\pm \frac {q^m+q^{-m}}{q-q^{-1}} \,
|m\rangle '',
$$
and
$$
Q^\pm _1(I_2)|0\rangle =\frac 1{q-q^{-1}} \, |1\rangle ',\ \ \
Q^\pm _1(I_2)|1\rangle ''=\frac 1{q-q^{-1}} \, |2\rangle '', $$
$$
Q^\pm _1(I_2)|m\rangle '=\frac 1{q-q^{-1}} \, |m+1\rangle '+
\frac 1{q-q^{-1}} \, |m-1\rangle ',\ \ \ m>0, $$
$$
Q^\pm _1(I_2)|m\rangle ''=\frac 1{q-q^{-1}} \, |m+1\rangle ''+
\frac 1{q-q^{-1}} \, |m-1\rangle '',\ \ \ m>1.
$$
Thus, the subspaces $V_1$ and $V_2$ are invariant with respect
to the representation $Q^+_1$ (and the representation $Q^-_1$).
We denote the subrepresentations
of $Q^\pm _1$ realized on $V_1$ and $V_2$ by $Q^{1,\pm}_1$ and
$Q^{2,\pm}_1$, respectively.

The eigenvalues of the operators $Q^\pm _{\sqrt q}(I_1)$ are
$$
\cdots ,\ q^{-3/2}+q^{3/2},\ q^{-1/2}+q^{1/2}, \ q^{1/2}+q^{-1/2},\
q^{3/2}+q^{-3/2},\ \cdots .
$$
Thus, every spectral point has multiplicity 2.
We denote by $W_1$ and $W_2$ the vector subspaces of ${\cal H}$ spanned
by the basis vectors
$$\textstyle
|1/2\rangle '=|0\rangle -|-1\rangle  , \
|3/2\rangle '=|1\rangle -|-2\rangle , \cdots ,
|m+\frac 12 \rangle '=|m \rangle -|-m-1 \rangle ,\cdots
$$
and
$$\textstyle
|1/2\rangle ''=|0\rangle +|-1\rangle , \ |3/2\rangle ''=|1\rangle +
|-2\rangle , \cdots , |m+\frac 12 \rangle ''=
|m \rangle +|-m-1\rangle ,\cdots ,
$$
respectively. These basis vectors are eigenvectors of the operator
$Q^\pm _{\sqrt q}(I_1)$:
$$
Q^\pm _{\sqrt q}(I_1)| m+{\textstyle \frac 12} \rangle ' =\pm
\frac {q^{m+1/2}+q^{-m-1/2}}
{q-q^{-1}} \, |m+{\textstyle \frac 12} \rangle ', $$
$$
Q^\pm _{\sqrt q}(I_1)| m+{\textstyle \frac 12} \rangle '' =\pm
\frac {q^{m+1/2}+q^{-m-1/2}}
{q-q^{-1}} \, |m+{\textstyle \frac 12} \rangle ''
$$
and
$$
Q^\pm _{\sqrt q}(I_2) |{\textstyle \frac 12}\rangle '=
- \frac 1{q-q^{-1}}\,  |{\textstyle \frac 12}\rangle '+
 \frac 1{q-q^{-1}}\,  |{\textstyle \frac 32}\rangle ',$$
$$
Q^\pm _{\sqrt q}(I_2) |m+{\textstyle \frac 12}\rangle '=
 \frac 1{q-q^{-1}}\,  |m+{\textstyle \frac 32}\rangle '+
\frac 1{q-q^{-1}}\,  |m-{\textstyle \frac 12}\rangle ',\ \ \ m>0,$$
$$
Q^\pm _{\sqrt q}(I_2) |{\textstyle \frac 12}\rangle ''=
 \frac 1{q-q^{-1}}\,  |{\textstyle \frac 12}\rangle ''+
 \frac 1{q-q^{-1}}\,  |{\textstyle \frac 32}\rangle '',$$
$$
Q^\pm _{\sqrt q}(I_2) |m+{\textstyle \frac 12}\rangle ''=
 \frac 1{q-q^{-1}}\,  |m+{\textstyle \frac 32}\rangle ''+
\frac 1{q-q^{-1}}\,  |m-{\textstyle \frac 12}\rangle '',\ \ \ m>0.
$$
Thus, the subspaces $W_1$ and $W_2$ are invariant with respect to the
representations $Q^\pm _{\sqrt q}$. We denote the subrepresentations of
$Q^\pm _{\sqrt q}$ realized on $W_1$ and $W_2$ by
$Q_{\sqrt q}^{1,\pm}$ and $Q_{\sqrt q}^{2,\pm}$, respectively.
Proposition is proved.
\medskip

{\bf Theorem 3.} {\it The representations
$Q^{1,\pm}_1$, $Q^{2,\pm}_1$, $Q_{\sqrt q}^{1,\pm}$ and
$Q_{\sqrt q}^{2,\pm}$ are irreducible and pairwise
nonequivalent. For any admissible value of $\lambda$
the representation $Q^+_\lambda$
(as well as the representation $Q^-_\lambda$) is
not equivalent to some of these representations.}
\medskip

{\sl Proof.} Proof is similar to that of Theorem 1 if to take into account
spectra of the operators $Q^{1,\pm}_1(I_1)$,
$Q^{2,\pm}_1(I_1)$, $Q_{\sqrt q}^{1,\pm}(I_1)$,
$Q_{\sqrt q}^{2,\pm}(I_1)$, $Q^\pm _\lambda (I_1)$ and traces of
the operators $Q^{1,\pm}_1(I_2)$, $Q^{2,\pm}_1(I_2)$,
$Q_{\sqrt q}^{1,\pm}(I_2)$, $Q_{\sqrt q}^{2,\pm}(I_2)$.
\medskip

   It will be proved in a separate paper that every irreducible infinite
dimensional representation of $U_q({\rm so}_3)$ is equivalent to one
of the representations described in this and previous sections.
\bigskip

\noindent
{\sf VIII. FINITE DIMENSIONAL REPRESENTATIONS OF
${\hat U}_q({\rm sl}_2)$: $q$ IS A ROOT OF UNITY}
\medskip

\noindent
Everywhere below $q$ is a root of unity, that is, there is a smallest
positive integer $p$ such that $q^p=1$. We suppose that $p\ne 1,\, 2$.
We introduce the number $p'$ setting $p'=p$ if $p$
is odd and $p'=p/2$ if $p$ is even.

As in the case of the algebra $U_q({\rm sl}_2)$ (see [11], Chapter 3), if
$q$ is a root of unity, then
$U_q({\rm so}_3)$ is a finite dimensional vector space
over the center of $U_q({\rm so}_3)$. If $q$ is a primitive root of unity,
then this assertion is stated in [5]. If $q$ is any root of unity, then
this assertion may be proved in the following way.
If $q^p=1$, then the center ${\cal C}$ of $U_q({\rm so}_3)$
contains the elements
$$
P_p=I_j^p+aI^{p-2}_j+bI^{p-4}_j+\cdots +dI_j^r ,\ \ \ j=1,2,3,
$$
where $r=0$ if $p$ is even and $r=1$ if $p$ is odd and $a,b,\cdots ,d$ are
certain fixed complex numbers expressed in terms of $q$.
(They are the polynomials $P$ defined in [5] if $q$ is a primitive root
of unity. Unfortunately, we could not find the explicit expressions for the
coefficients $a,b,\cdots ,d$. But note that $P_3=I_j^3+I_j$,
$P_4=I_j^4+I_j^2$ and $P_5=I_j^5+(1+(q+q^{-1})^{-1}I_j^3+(q+q^{-1})^{-1}
I_j$.) Therefore, $I_j^s$,
$s>n$, can be reduced to the linear combination of $I^i_j$, $i<n$, with
coefficients from the center ${\cal C}$. Now our assertion follows from
this and from Poincar\'e--Birkhoff--Witt theorem for $U_q({\rm so}_3)$.
\medskip

{\bf Theorem 4.} {\it If $q$ is a root of unity, then any irreducible
representation of $U_q({\rm so}_3)$ is finite dimensional.}
\medskip

{\sl Proof.} Let $T$ be an irreducible representation of $U_q({\rm so}_3)$.
Then $T$ maps central elements into scalar operators. Since the linear
space $U_q({\rm so}_3)$ is finite dimensional over the center ${\cal C}$
with the basis $I_1^kI_2^mI_3^n $, $k,m,n<p$, then
for any $a\in U_q({\rm so}_3)$ we have $T(a)=\sum _{k,m,n<p}
T(I_1^kI_2^mI_3^n)$. Hence, if ${\bf v}$ is a nonzero vector of the
representation space ${\cal V}$, then $T(U_q({\rm so}_3)){\bf v}={\cal V}$
and ${\cal V}$ is finite dimensional. Theorem is proved.
\medskip

Taking into account Theorem 4, below we consider only finite
dimensional representations of $U_q({\rm so}_3)$.

In order to find irreducible representations of $U_q({\rm so}_3)$ for
$q$ a root of unity, we use the same method as before, that is, we apply the
homomorphism $\psi$ from Proposition 2 and irreducible representations of
the algebra ${\hat U}_q({\rm sl}_2)$ for $q$ a root of unity.

Let us find irreducible representations of ${\hat U}_q({\rm sl}_2)$
for $q$ a root of unity. The quantum algebra $U_q({\rm sl}_2)$ for $q$ a root of
unity has the following irreducible representations (see [11], Subsection
3.3.2):
\medskip

(a) The representations $T^{(1)}_l$, $T^{(-1)}_l$, $T^{({\rm i})}_l$,
$T^{(-{\rm i})}_l$, $2l<p'$, given by the formulas (15)--(19).
\smallskip

(b) The representations $T_{ab\lambda}$, $a,b,\lambda \in {\bf C}$, $\lambda
\ne 0$, acting on a $p'$-dimensional vector space ${\cal H}$ with the basis
$| j\rangle$, $j=0,1,2,\cdots ,p'-1$, and given by the formulas
$$
T_{ab\lambda} (q^H)| i\rangle =q^{-i}\lambda | i\rangle ,\ \ \
T_{ab\lambda} (F)| p'-1\rangle =b | 0\rangle , \eqno (36) $$
$$
T_{ab\lambda} (F)| i\rangle = | i+1 \rangle , \ i<p'-1,\ \ \
T_{ab\lambda} (E)| 0\rangle =a | p'-1 \rangle , \eqno (37) $$
$$
T_{ab\lambda} (E)| i\rangle =\left( ab+[i]\, \frac {\lambda ^2q^{1-i}-
\lambda ^{-2}q^{i-1}}{q-q^{-1}} \right) \, | i-1\rangle ,\ \ i>0 .
\eqno (38)
$$
The representations $T_{ab\lambda}$ with $(a,b)=(0,0)$ and $\lambda =\pm q^n$,
$n=0,1,2,\cdots ,p'-2$, are reducible and must be taken out from this set.
\smallskip

(c) The representations $T'_{0b\lambda}$, $b,\lambda \in {\bf C}$, $\lambda
\ne 0$, acting on a $p'$-dimensional vector space ${\cal H}$ with the basis
$| j\rangle$, $j=0,1,2,\cdots ,p'-1$, and given by the formulas
$$
T'_{0b\lambda} (q^H)| i\rangle =q^{i}\lambda ^{-1}| i\rangle ,\ \ \
T'_{0b\lambda} (E)| p'-1\rangle =b | 0\rangle , \eqno (39) $$
$$
T'_{0b\lambda} (E)| i\rangle = | i+1 \rangle , \ i<p'-1,\ \ \
T_{0b\lambda} (F)| 0\rangle =0 , \eqno (40) $$
$$
T'_{0b\lambda} (F)| i\rangle =[i]\, \frac {\lambda ^2q^{1-i}-
\lambda ^{-2}q^{i-1}}{q-q^{-1}} \, | i-1\rangle ,\ \ i>0 .
\eqno (41)
$$
The representations $T'_{00\lambda}$ with $\lambda =\pm q^n$,
$n=0,1,2,\cdots ,p'-2$, are reducible and must be taken out from this set.
\medskip

\noindent
{\sl Remark 1.} In the set of representations (a)--(c) there exist
equivalent representations (see, for example, Propositions 3.17 and 3.18
in [11]).
\medskip

\noindent
{\sl Remark 2.} In [11], Subsection 3.3.2, irreducible representations of
the algebra generated by the elements $E, F, K:=q^{2H}, K^{-1}:=q^{-2H}\in
U_q({\rm sl}_2)$ are given. Clearly, this algebra is a subalgebra in
$U_q({\rm sl}_2)$. It is easy to generalize the results of
Subsection 3.3.2
in [11] for $U_q({\rm sl}_2)$. Let us note that the algebra $U_q({\rm sl}_2)$
has a unique automorphism $\varphi$ such that $\varphi (q^H)={\rm i}q^H$,
$\varphi (E)=-E$ and $\varphi (F)=F$. (If $q$ is not a root of unity, then
this automorphism transforms the representations $T^{(1)}_l$ to the
representations $T^{({\rm i})}_l$, respectively.)
Therefore, the mapping
${\tilde T}_{-a,b,\lambda}=T_{ab\lambda}\circ \varphi$ is also a
representation of $U_q({\rm sl}_2)$. We have
$$
{\tilde T}_{ab\lambda} (q^H)| i\rangle ={\rm i}q^{-i}\lambda |i\rangle ,\ \ \
{\tilde T}_{ab\lambda} (F)| p'-1\rangle =b | 0\rangle , \eqno (42) $$
$$
{\tilde T}_{ab\lambda} (F)| i\rangle = | i+1 \rangle , \ i<p'-1,\ \ \
{\tilde T}_{ab\lambda} (E)| 0\rangle =a | p'-1 \rangle , \eqno (43) $$
$$
{\tilde T}_{ab\lambda} (E)| i\rangle =\left( ab-[i]\, \frac
{\lambda ^2q^{1-i}-
\lambda ^{-2}q^{i-1}}{q-q^{-1}} \right) \, | i-1\rangle ,\ \ i>0 .
\eqno (44)
$$
However, it is easy to see by comparing (36)--(38) with (42)--(44) that
the representation ${\tilde T}_{ab\lambda}$ is equivalent to
$T_{a,b,{\rm i}\lambda}$. This means that for $q$ a root of unity
we do not obtain
new representations of $U_q({\rm sl}_2)$ from $T_{ab\lambda}$ applying
the automorphism
$\varphi$ as in the case of the representations $T^{(1)}_l$.
\medskip

 We have described irreducible representations of the algebra
$U_q({\rm sl}_2)$. Now we wish to extend these representations to
obtain representations of
the algebra ${\hat U}_q({\rm sl}_2)$ by using the relation
$$
T((q^kq^H+q^{-k}q^{-H})^{-1}):=(q^kT(q^H)+q^{-k}T(q^{-H}))^{-1}.
$$
Clearly, only those irreducible representations $T$ of $U_q({\rm sl}_2)$
can be extended to ${\hat U}_q({\rm sl}_2)$ for which the operators
$q^kT(q^H)+q^{-k}T(q^{-H})$ are invertible. From formulas (15)--(19) it is
clear that these operators are always invertible for the
irreducible representations
$T_l^{(1)}$, $T_l^{(-1)}$, $l=0,\frac 12 ,1,\frac 32 ,\cdots ,\frac {p'-1}2$,
and for the irreducible representations
$T_l^{({\rm i})}$, $T_l^{(-{\rm i})}$, $l=\frac 12 ,\frac 32 ,
\frac 52 ,\cdots ,\frac {p'-1}2$ (or $\frac {p'-2}2$). (For the
representations $T_l^{({\rm i})}$, $T_l^{(-{\rm i})}$, $l=0,1,2,
\cdots $, some of these operators are not invertible since they have
zero eigenvalue.) We denote the extended representations by the same symbols
$T_l^{(1)}$, $T_l^{(-1)}$, $T_l^{({\rm i})}$, $T_l^{(-{\rm i})}$,
respectively.

Similarly, the representation $T_{ab\lambda}$ (and the representation
$T'_{0b\lambda}$) can be extended to a representation of the algebra
${\hat U}_q({\rm sl}_2)$ if and only if $\lambda \ne \pm {\rm i}q^k$,
$k\in {\bf Z}$.
\medskip

{\bf Proposition 7.} {\it The algebra ${\hat U}_q({\rm sl}_2)$
for $q$ a root of unity has the irreducible representations
$T_l^{(1)}$, $T_l^{(-1)}$, $l=0,\frac 12 ,1,\frac 32 ,\cdots ,\frac {p'-1}2$,
the irreducible representations
$T_l^{({\rm i})}$, $T_l^{(-{\rm i})}$, $l=\frac 12 ,\frac 32 ,
\frac 52 ,\cdots \frac {p'-1}2$ (or $\frac {p'-2}2$), and the irreducible
representations $T_{ab\lambda}$, $T'_{0b\lambda}$, $\lambda \ne \pm {\rm i}
q^k$, $k\in {\bf Z}$. Any irreducible representation of
${\hat U}_q({\rm sl}_2)$ for $q$ a root of unity is equivalent to one of
these representations.}
\bigskip

\noindent
{\sf IX. REPRESENTATIONS OF
$U_q({\rm so}_3)$ FOR $q$ A ROOT OF UNITY
OBTAINED FROM THOSE OF ${\hat U}_q({\rm sl}_2)$}
\medskip

As in Section IV, we shall obtain representations of $U_q({\rm so}_3)$
for $q$ a root of unity by applying the homomorphism $\psi$ from
Proposition 2. Namely, if $T$ is a representation of
${\hat U}_q({\rm sl}_2)$,
then
$$
R=T\circ \psi \eqno (45)
$$
is a representation of $U_q({\rm so}_3)$. As in Section IV, application of
this method to the pair of the irreducible representations $T_l^{(1)}$ and
$T_l^{(-1)}$ of ${\hat U}_q({\rm sl}_2)$ leads to the same representation of
$U_q({\rm so}_3)$ which will be denoted by $R_l^{(1)}$. Applying the
formula (45) to the irreducible representations $T_l^{({\rm i})}$ and
$T_l^{(-{\rm i})}$ of ${\hat U}_q({\rm sl}_2)$ give the representations
of $U_q({\rm so}_3)$ which will be denoted by $R_l^{({\rm i})}$ and
$R_l^{(-{\rm i})}$, respectively.
\medskip

{\bf Proposition 8.} {\it The representations $R^{(1)}_l$ of $U_q({\rm so}_3)$
are irreducible. The representations $R_l^{({\rm i})}$ and $R_l^{(-{\rm i})}$
are reducible.}
\medskip

{\sl Proof} of this proposition is the same as that of Proposition 3.
\medskip

Repeating word-by-word the reasoning of Section IV, we decompose the
representations $R_l^{({\rm i})}$ and $R_l^{(-{\rm i})}$ into the direct
sums of representations of $U_q({\rm so}_3)$ which are denoted by
$R_n^{(\pm {\rm i},+)}$ and $R_n^{(\pm {\rm i},-)}$:
$$
R_l^{({\rm i})}=R_n^{({\rm i},+)}\oplus R_n^{({\rm i},-)},\ \ \
R_l^{(-{\rm i})}=R_n^{(-{\rm i},+)}\oplus R_n^{(-{\rm i},-)},\ \ \
n=l+\frac 12 .
$$
Moreover, the representations
$R_n^{(\pm {\rm i},+)}$ and $R_n^{(\pm {\rm i},-)}$ are given in the
appropriate bases $|1\rangle$, $|2\rangle$, $\cdots $, $|n\rangle$ by the
corresponding formulas of Section IV.
\medskip

{\bf Theorem 5.} {\it The representations
$R_n^{({\rm i},+)}$, $R_n^{({\rm i},-)}$,
$R_n^{(-{\rm i},+)}$, $R_n^{(-{\rm i},-)}$, $n=1,2,3,\cdots ,\frac {p'}2$
(or $\frac {p'-1}2$) are irreducible and pairwise nonequivalent. For any
$l$, $l=0,\frac 12 ,1,\frac 32 ,\cdots ,\frac {p'-1}2$, the representation
$R_l^{(1)}$ is not equivalent to some of these representations.}
\medskip

{\sl Proof} is the same as that of Theorem 1.
\medskip

Now we apply formula (45) to the representations $T_{ab\lambda}$ and
$T'_{0b\lambda}$. As a result, we obtain the representations
$$
R_{ab\lambda}=T_{-a,b,-{\rm i}\lambda} \circ \psi ,\ \ \
R'_{0b\lambda}=T'_{0,b,-{\rm i}\lambda}
$$
given in the bases $|j\rangle$, $j=0,1,2,\cdots ,p'-1$, by the formulas
$$
R_{ab\lambda}(I_1)|i\rangle =\frac {-1}{q-q^{-1}} (q^{-i}\lambda +
q^i\lambda ^{-1})|i\rangle , \eqno (46) $$
$$
R_{ab\lambda}(I_2)|0\rangle =\frac {{\rm i}}
{\lambda -\lambda ^{-1}} \left( a| p'-1\rangle +
|1\rangle \right) , \eqno (47) $$
$$
R_{ab\lambda}(I_2)|p'-1\rangle =\frac {{\rm i}}
{q^{-p'+1}\lambda -q^{p'-1}\lambda ^{-1}} \Biggl\{ b| 0\rangle +
\qquad\qquad\qquad$$
$$\qquad\qquad\qquad
+\left( ab+[p'-1]\frac {q^{-p'+2}\lambda ^2-q^{p'-2}\lambda ^{-2}}{q-q^{-1}}
\right) |p'-2\rangle  \Biggr\} , \eqno (48) $$
$$
R_{ab\lambda}(I_2)|i \rangle =\frac {{\rm i}}
{q^{-i}\lambda -q^i\lambda ^{-1}} \Biggl\{ \left( ab+[i]
\frac {q^{-i+1}\lambda ^2-q^{i-1}\lambda ^{-2}}{q-q^{-1}}\right)
| i-1\rangle + $$
$$\qquad\qquad\qquad\qquad\qquad\qquad
+|i+1\rangle  \Biggr\} ,\ \ \  0<i<p'-1 . \eqno (49)
$$
and by the formulas
$$
R'_{0b\lambda}(I_1)|i\rangle =\frac {1}{q-q^{-1}} (q^{-i}\lambda +
q^i\lambda ^{-1})|i\rangle ,  \ \ \
R'_{0b\lambda}(I_2)|0\rangle =\frac {-{\rm i}}
{\lambda -\lambda ^{-1}} |1\rangle ,  $$
$$
R'_{0b\lambda}(I_2)|p'-1\rangle =\frac {-{\rm i}}
{q^{-p'+1}\lambda -q^{p'-1}\lambda ^{-1}} \Biggl( b| 0\rangle
+[p'-1]\frac {q^{-p'+2}\lambda ^2-q^{p'-2}\lambda ^{-2}}{q-q^{-1}}
|p'-2\rangle  \Biggr) ,  $$
$$
R'_{0b\lambda}(I_2)|i \rangle =\frac {-{\rm i}}
{q^{-i}\lambda -q^i\lambda ^{-1}} \Biggl( | i+1\rangle +
[i] \frac {q^{-i+1}\lambda ^2-q^{i-1}\lambda ^{-2}}{q-q^{-1}}
| i-1\rangle \Biggr) , $$
$$\qquad\qquad\qquad\qquad\qquad\qquad\qquad\qquad\qquad
0<i<p'-1 .
$$
The operators $R_{ab\lambda}(I_3)$ and $R'_{0b\lambda}(I_3)$ can be
calculated by means of the relation
$$
R(I_3)=q^{1/2}R(I_1)R(I_2)-q^{-1/2}R(I_2)R(I_1).
$$
Recall that the representations $R_{ab\lambda}$ and $R'_{0b\lambda}$ are
determined for $\lambda \ne 0$ and $\lambda \ne \pm q^k$, $k\in {\bf Z}$.

It is seen from the above formulas that
$$
R'_{0b\lambda}(I_1)= R_{0,b,-\lambda}(I_1),\ \ \
R'_{0b\lambda}(I_2)= R_{0,b,-\lambda}(I_2),
$$
that is, the
representations $R_{0,b,-\lambda}$ and $R'_{0b\lambda}$ are equivalent. For
this reason, we consider below only the representations $R_{ab\lambda}$.

In order to study the representations $R_{ab\lambda}$ of $U_q({\rm so}_3)$
we consider the spectrum of the operator $R_{ab\lambda}(I_1)$. It
coincides with the set of points
$$
-\frac {\lambda +\lambda ^{-1}}{q-q^{-1}},\ \
-\frac {q^{-1}\lambda +q\lambda ^{-1}}{q-q^{-1}},\ \
-\frac {q^{-2}\lambda +q^2\lambda ^{-1}}{q-q^{-1}},\ \cdots \ ,
-\frac {q^{1-p'}\lambda +q^{p'-1}\lambda ^{-1}}{q-q^{-1}}.
\eqno (50)
$$
It is easy to see that there exist coinciding points in this set if and
only if $\lambda$ is equal to one of the numbers
$$
\pm q^{1/2},\ \  \pm q^{3/2},\ \ \pm q^{5/2},\ \cdots \ ,
\pm q^{(p'-1)/2} \ ({\rm or}\ \pm q^{(p'-2)/2}).
$$
(Here we have to take $\pm q^{(p'-1)/2}$ if $p'$ is even and
$\pm q^{(p'-2)/2}$ if
$p'$ is odd.) Moreover, the set (50) splits into pairs of coinciding points
if and only if $\lambda =\pm q^{(p'-1)/2}$. In all other cases there exists
at least one spectral point which coincides with no other point.
In particular, if $\lambda =\pm q^{(p'-2)/2}$, then in this set there exists
only one eigenvalue with multiplicity 1. In all other cases there are more
than one eigenvalues with multiplicity 1.
\medskip

{\bf Proposition 9.} {\it If $\lambda \ne \pm q^{(p'-1)/2}$ for even $p'$ and
$\lambda \ne \pm q^{(p'-2)/2}$ for odd $p'$, then the
representation $R_{ab\lambda}$ is irreducible.}
\medskip

{\sl Proof.} Let $\lambda \ne \pm q^{(p'-1)/2}$ for even $p'$ and
$\lambda \ne \pm q^{(p'-2)/2}$ for odd $p'$.
We distinguish two cases:
when the spectrum of the operator $R_{ab\lambda}(I_1)$ is simple and when
there exists at list one spectral point of this operator having multiplicity
2. In the first case the proof is the same as the first part of the proof of
Proposition 3. For the second case, we give a proof only for $\lambda =
q^{1/2}$. (Proofs for other values of $q$ are similar.) Then in the set
(50) there are only two coinciding points
$-\frac {\lambda +\lambda ^{-1}}{q-q^{-1}}$ and
$-\frac {q^{-1}\lambda +q\lambda ^{-1}}{q-q^{-1}}$ corresponding
to the eigenvectors $|0\rangle$ and $|1\rangle$. Let $V$ be an invariant
subspace of the representation space ${\cal H}$. As in the proof of
Proposition 3, it is shown that $V$ is a linear span of eigenvectors
of the operator $R_{ab\lambda}(I_1)$, that is, a certain part of the
vectors $|i\rangle$, $i\ne 0,\, 1$, $\alpha _0|0\rangle +\alpha _1|1\rangle$,
$\beta _0|0\rangle +\beta _1|1\rangle$ constitutes a basis of $V$. Let $V$
contain some basis vector $|j\rangle$. Then as in the proof of Proposition
3, acting successively upon $|j\rangle$ by certain linear combinations
of the operators $R_{ab\lambda}(I_2)$ and $R_{ab\lambda}(I_3)$ we generate
all the vectors $|i\rangle$, $i=0,1,\cdots ,\frac 12 (p'-1)$. This means
that $V={\cal H}$ and the representation $R_{ab\lambda}$ is irreducible.
If $V$ contains no vector $|j\rangle$, $j\ne 0,\, 1$, then some
linear combination $\alpha _0|0\rangle +\alpha _1|1\rangle$ belongs
to $V$. Then the vector ${\bf v}=R_{ab\lambda}(I_2)
(\alpha _0|0\rangle +\alpha _1|1\rangle)$ belongs to $V$. Since ${\bf v}$
contains the summand $\alpha |2\rangle$ with nonzero coefficient $\alpha$,
then $|2\rangle \in V$. This is a contradiction. Hence, the representation
$R_{ab\lambda}$ is irreducible. Proposition is proved.
\medskip

Let $p'$ be even. Let us study the representations $R_{ab\lambda}$ for
$\lambda =\pm q^{(p'-1)/2}$. For $\lambda =q^{(p'-1)/2}$ we have
$$
R_{ab\lambda}(I_1)|i\rangle =\frac {-1}{q-q^{-1}} (q^{-i+(p'-1)/2}+
q^{i-(p'-1)/2})\, |i\rangle , \eqno (51) $$
$$
R_{ab\lambda}(I_2)|0\rangle =c_{(p'-1)/2}\left( a|p'-1\rangle +|1\rangle
\right) , \eqno (52)  $$
$$
R_{ab\lambda}(I_2)|p'-1\rangle =-c_{(p'-1)/2}\left( (ab+[p'-1]^2)
|p'-2\rangle +b|0\rangle \right) , \eqno (53)  $$
$$
R_{ab\lambda}(I_2)|i\rangle =c_{-i+(p'-1)/2}\left( (ab+[i]^2)
|i-1\rangle +|i+1\rangle \right) , \eqno (54)
$$
where
$$
c_j=\frac {\rm i}{q^j-q^{-j}} .
$$
The operator $R_{a,b,(p'-1)/2}(I_1)$ has the spectrum
$$
\frac {-1}{q-q^{-1}} (q^{-i+(p'-1)/2}+
q^{i-(p'-1)/2}), \ \ \ i=0,1,2,\cdots ,p'-1 ,
$$
that is, if $p'$ is even, then all spectral points are of multiplicity 2.

We assume that $ab\ne -[j]^2$, $j=0,1,\cdots ,p'-1$, and
go over from the basis $\{ |i\rangle\}$ to the basis $\{ |i\rangle
^\circ \}$, where
$$
| i\rangle ^\circ =\prod _{j=0}^i (ab+[j]^2)^{-1/2} | i\rangle ,\ \ \
i=0,1,2,\cdots ,p'-1 .
$$
Then the formula (51) does not change and the formulas (52)--(54) turn into
$$
R_{ab\lambda}(I_2)|0\rangle ^\circ =c_{(p'-1)/2}
\Biggl( a\prod _{j=1}^{p'-1}(ab+[j]^2)^{1/2}\, |p'-1\rangle ^\circ +
(ab+1)^{1/2} |1\rangle ^\circ \Biggr) ,   $$
$$
R_{ab\lambda}(I_2)|p'-1\rangle ^\circ =-c_{(p'-1)/2}\Biggl( (ab+1)^{1/2}
|p'-2\rangle ^\circ +\frac b{\prod _{j=1}^{p'-1}(ab+[j]^2)^{1/2}}
|0\rangle ^\circ \Biggr) ,   $$
$$
R_{ab\lambda}(I_2)|i\rangle ^\circ =c_{-i+(p'-1)/2}\left( (ab+[i]^2)^{1/2}
|i-1\rangle ^\circ +(ab+[i+1]^2)^{1/2} |i+1\rangle ^\circ \right) .
$$
We split the representation
space ${\cal H}$ into the direct sum of two linear subspaces
${\cal H}_1$ and ${\cal H}_2$ spanned by the basis vectors $|j\rangle '$,
$j=0,1,2,\cdots ,\frac 12 (p'-2)$, and $|j\rangle ''$,
$j=0,1,2,\cdots ,\frac 12 (p'-2)$, where
$$
|j\rangle '=|j\rangle ^\circ +{\rm i}(-1)^{-j-1+p'/2}|p'-j-1\rangle ^\circ ,
\ \ \
|j\rangle ''=|j\rangle ^\circ +{\rm i}(-1)^{-j+p'/2}|p'-j-1\rangle ^\circ .
$$
Then as in Section IV, we derive
$$
R_{a,b,(p'-1)/2}(I_1)|j\rangle '=\frac {-1}{q-q^{-1}}
(q^{-j+(p'-1)/2}+q^{j-(p'-1)/2})\, |j\rangle ', $$
$$
R_{a,b,(p'-1)/2}(I_1)|j\rangle ''=\frac {-1}{q-q^{-1}}
(q^{-j+(p'-1)/2}+q^{j-(p'-1)/2})\, |j\rangle ''
$$
for the operator $R_{a,b,(p'-1)/2}(I_1)$ and
$$
R_{a,b,(p'-1)/2}(I_2)|j\rangle '=c_{-j+(p'-1)/2}\left( (ab+[j+1]^2)^{1/2}
|j+1\rangle '+(ab+[j]^2)^{1/2} |j-1\rangle ' \right) ,$$
$$
R_{a,b,(p'-1)/2}(I_2)|j\rangle ''=c_{-j+(p'-1)/2}\left( (ab+[j+1]^2)^{1/2}
|j+1\rangle ''+(ab+[j]^2)^{1/2} |j-1\rangle '' \right) ,
$$
where $j\ne 0,\frac {p'}2 -1$,
$$
R_{a,b,(p'-1)/2}(I_2)|{\textstyle \frac {p'}2} -1\rangle '=
\frac 1{q^{1/2}-q^{-1/2}} (ab+[{\textstyle \frac {p'}2}]^2)^{1/2}
|{\textstyle \frac {p'}2} -1\rangle ' +$$
$$\qquad\qquad\qquad\qquad\qquad\qquad
+\frac {\rm i}{q^{1/2}-q^{-1/2}}
(ab+[{\textstyle \frac {p'}2} -1]^2)^{1/2}
|{\textstyle \frac {p'}2} -2\rangle ' ,$$
$$
R_{a,b,(p'-1)/2}(I_2)|{\textstyle \frac {p'}2} -1\rangle ''=-
\frac 1{q^{1/2}-q^{-1/2}} (ab+[{\textstyle \frac {p'}2}]^2)^{1/2}
|{\textstyle \frac {p'}2} -1\rangle '' +$$
$$\qquad\qquad\qquad\qquad\qquad\qquad
+\frac {\rm i}{q^{1/2}-q^{-1/2}}
(ab+[{\textstyle \frac {p'}2} -1]^2)^{1/2}
|{\textstyle \frac {p'}2} -2\rangle '' ,$$
$$
R_{a,b,(p'-1)/2}(I_2)|0\rangle '=c_{(p'-1)/2}\Biggl(
a\prod _{j=1}^{p'-1} (ab+[j]^2)^{1/2} \, |p'-1\rangle ^\circ
+(ab+1)^{1/2}|1\rangle ^\circ \Biggr) $$
$$
-{\rm i}(-1)^{(p'-2)/2}c_{(p'-1)/2} \Biggl( (ab+1)^{1/2} |p'-2\rangle ^\circ
+\frac b{\prod _{j=1}^{p'-1} (ab+[j]^2)^{1/2}}\, |0\rangle ^\circ \Biggr) .
$$
When
$$
a\prod _{j=1}^{p'-1} (ab+[j]^2)^{1/2}=\frac b
{\prod _{j=1}^{p'-1} (ab+[j]^2)^{1/2}} , \eqno (55)
$$
then the last relation reduces to
$$
R_{a,b,(p'-1)/2}(I_2)|0\rangle '=\frac {(-1)^{(p'-2)/2}}{q^{(p'-1)/2}-
q^{-(p'-1)/2}}\, a\prod _{j=1}^{p'-1} (ab+[j]^2)^{1/2} \, |0\rangle '+
$$
$$\qquad\qquad\qquad\qquad\qquad
+c_{(p'-1)/2}(ab+1)^{1/2}|1\rangle '.
$$
Similarly, if the condition (55) is fulfilled, then
$$
R_{a,b,(p'-1)/2}(I_2)|0\rangle ''=\frac {(-1)^{p'/2}}{q^{(p'-1)/2}-
q^{-(p'-1)/2}}\, a\prod _{j=1}^{p'-1} (ab+[j]^2)^{1/2} \, |0\rangle ''+
$$
$$\qquad\qquad\qquad\qquad\qquad
+c_{(p'-1)/2}(ab+1)^{1/2}|1\rangle ''.
$$
Thus, the subspaces ${\cal H}_1$ and ${\cal H}_2$ are invariant with respect
to the representation $R_{a,b,(p'-1)/2}$ if the condition (55) is fulfilled.
We denote the corresponding subrepresentations by $R^{1,+}_{a,b,(p'-1)/2}$
and $R^{2,+}_{a,b,(p'-1)/2}$, respectively.

Similarly, if $\lambda =-q^{(p'-1)/2}$, then
$$
R_{a,b,-(p'-1)/2}(I_1)=- R_{a,b,(p'-1)/2}(I_1) ,\ \ \
R_{a,b,-(p'-1)/2}(I_2)=- R_{a,b,(p'-1)/2}(I_2)
$$
and the subspaces
${\cal H}_1$ and ${\cal H}_2$ are invariant with respect
to the representation $R_{a,b,-(p'-1)/2}$ if the condition (55) is fulfilled.
We denote the corresponding subrepresentations by $R^{1,-}_{a,b,-(p'-1)/2}$
and $R^{2,-}_{a,b,-(p'-1)/2}$, respectively.
\medskip

{\bf Proposition 10.} {\it Let the condition (55) is satisfied. Then the
representations $R^{i,+ }_{a,b,(p'-1)/2}$ and
$R^{i,-}_{a,b,-(p'-1)/2}$, $i=1,2$, of the algebra $U_q({\rm so}_3)$
are irreducible and pairwise nonequivalent. If the condition (55) is not
satisfied, then the representations
$R_{a,b,(p'-1)/2}$ and $R_{a,b,-(p'-1)/2}$ are irreducible.}
\medskip

{\sl Proof} is similar to that of the previous propositions and we omit it.
\medskip

Remark that the representations $R^{i,+ }_{a,b,(p'-1)/2}$ and
$R^{i,- }_{a,b,-(p'-1)/2}$, $i=1,2$, have two nonzero diagonal matrix
elements $\langle \frac {p'}2 -1 | R|\frac {p'}2 -1 \rangle$ and
$\langle 0 |R|0\rangle$.

Let now $p'$ be odd and $\lambda =q^{(p'-2)/2}$. For this value of
$\lambda$ we have
$$
R_{ab\lambda}(I_1)|i\rangle =\frac {-1}{q-q^{-1}} (q^{-i+(p'-2)/2}+
q^{i-(p'-2)/2})\, |i\rangle ,  $$
$$
R_{ab\lambda}(I_2)|0\rangle =c_{(p'-2)/2}\left( a|p'-1\rangle +|1\rangle
\right) ,   $$
$$
R_{ab\lambda}(I_2)|p'-1\rangle =-c_{p'/2}\left( (ab+\epsilon [p'-1][p'])
|p'-2\rangle +b|0\rangle \right) ,   $$
$$
R_{ab\lambda}(I_2)|i\rangle =c_{-i+(p'-2)/2}\left( (ab+\epsilon [i][i+1])
|i-1\rangle +|i+1\rangle \right) ,
$$
where $\epsilon =1$ for $p'=p/2$, $\epsilon =-1$ for $p'=p$ and $c_j$ is
such as in (51)--(54). The operator $R_{a,b,(p'-2)/2}(I_1)$ has the spectrum
$$
\frac {-1}{q-q^{-1}} (q^{-i+(p'-2)/2}+
q^{i-(p'-2)/2}), \ \ \ i=0,1,2,\cdots ,p'-1 ,
$$
that is, all spectral points are of multiplicity 2 except for the point
$-(q^{p'/2}+q^{-p'/2})/(q-q^{-1})$ which is of multiplicity 1.

We assume that $ab\ne -\epsilon [j][j+1]$, $j=0,1,\cdots ,p'-1$, and
go over from the basis $\{ |i\rangle\}$ to the basis $\{ |i\rangle
^\circ \}$, where
$$
| i\rangle ^\circ =\prod _{j=0}^i (ab+\epsilon [j][j+1])^{-1/2} | i\rangle ,
\ \ \ i=0,1,2,\cdots ,p'-1 .
$$
Then
$$
R_{ab\lambda}(I_1)|i\rangle ^\circ =\frac {-1}{q-q^{-1}} (q^{-i+(p'-2)/2}+
q^{i-(p'-2)/2})\, |i\rangle ^\circ ,  $$
$$
R_{ab\lambda}(I_2)|0\rangle ^\circ =c_{(p'-2)/2} \Biggl( a\prod _{j=1}^{p'-1}
(ab+\epsilon [j][j+1])^{1/2}\, |p'-1\rangle ^\circ
+ (ab+\epsilon [2])^{1/2} |1\rangle ^\circ \Biggr) ,   $$
$$
R_{ab\lambda}(I_2)|p'-1\rangle ^\circ =-c_{p'/2}( (ab+\epsilon
[p'-1][p'])^{1/2} |p'-2\rangle ^\circ +$$
$$\qquad\qquad\qquad\qquad\qquad
+ b\prod _{j=1}^{p'-1}(ab+\epsilon [j][j+1])^{-1/2}
|0\rangle ^\circ ) ,   $$
$$
R_{ab\lambda}(I_2)|i\rangle ^\circ =c_{-i+(p'-2)/2}( (ab+\epsilon
[i][i+1])^{1/2} |i-1\rangle ^\circ  +$$
$$\qquad\qquad\qquad\qquad\qquad
+(ab+\epsilon [i+1][i+2])^{1/2} |i+1\rangle ^\circ ) ,
$$
where $\lambda =q^{(p'-2)/2}$. Let ${\cal H}_1$ and ${\cal H}_2$
be two linear subspaces of the representation space ${\cal H}$
spanned by the basis vectors
$$
|j\rangle '=|j\rangle ^\circ +{\rm i}(-1)^j|p'-j-2\rangle ^\circ ,
\ \ \ j=0,1,2,\cdots ,\frac {p'-3}2 ,
$$
and the basis vectors
$$
|j\rangle ''=|j\rangle ^\circ +{\rm i}(-1)^{j+1}|p'-j-2\rangle ^\circ ,\ \ \
j=0,1,2,\cdots ,\frac {p'-3}2 ,
$$
respectively. Then the operator $R_{a,b,(p'-2)/2}(I_1)$ acts on the basis
elements $|j\rangle '$ and $|j\rangle ''$ as on the vectors $|j\rangle$ and
$$
R_{a,b,(p'-1)/2}(I_2)|j\rangle '=c_{-j+(p'-2)/2}\left( (ab+\epsilon
[j+1][j+2])^{1/2} |j+1\rangle ' \right.$$
$$\qquad\qquad\qquad\qquad\qquad   \left.
+(ab+\epsilon [j][j+1])^{1/2} |j-1\rangle '\right) ,$$
$$
R_{a,b,(p'-2)/2}(I_2)|j\rangle ''=c_{-j+(p'-2)/2}\left( (ab+\epsilon
[j+1][j+2])^{1/2} |j+1\rangle '' \right. $$
$$\qquad\qquad\qquad\qquad\qquad    \left.
+(ab+\epsilon [j][j+1])^{1/2} |j-1\rangle ''\right) ,
$$
where $j\ne 0,\frac {p'-3}2$,
$$
R_{a,b,(p'-2)/2}(I_2)|{\textstyle \frac {p'-3}2} \rangle '=
\frac {(-1)^{(p'-3)/2}}{q^{1/2}-q^{-1/2}}
(ab+\epsilon [{\textstyle \frac {p'-1}2}] [{\textstyle \frac {p'+1}2}])^{1/2}
|{\textstyle \frac {p'-3}2} \rangle ' +$$
$$\qquad\qquad\qquad\qquad
+\frac {\rm i}{q^{1/2}-q^{-1/2}}
(ab+\epsilon [{\textstyle \frac {p'-1}2}] [{\textstyle \frac {p'-3}2}])^{1/2}
|{\textstyle \frac {p'-5}2} \rangle ' ,$$
$$
R_{a,b,(p'-2)/2}(I_2)|{\textstyle \frac {p'-3}2} \rangle ''=-
\frac {(-1)^{(p'-3)/2}}{q^{1/2}-q^{-1/2}}
(ab+\epsilon [{\textstyle \frac {p'-1}2}] [{\textstyle \frac {p'+1}2}])^{1/2}
|{\textstyle \frac {p'-3}2} \rangle '' +$$
$$\qquad\qquad\qquad\qquad
+\frac {\rm i}{q^{1/2}-q^{-1/2}}
(ab+\epsilon [{\textstyle \frac {p'-1}2} ][{\textstyle \frac {p'-3}2} ])^{1/2}
|{\textstyle \frac {p'-5}2} \rangle '' ,$$
$$
R_{a,b,(p'-2)/2}(I_2)|0\rangle '=c_{(p'-2)/2}\Biggl(
a\prod _{j=1}^{p'-1} (ab+\epsilon [j][j+1])^{1/2} \, |p'-1\rangle ^\circ
+(ab+\epsilon [2])^{1/2}|1\rangle ^\circ  $$
$$
-{\rm i} (ab+\epsilon [2])^{1/2} |p'-1\rangle ^\circ -{\rm i}
(ab+\epsilon [p'-2][p'-1])^{1/2}\, |p'-3 \rangle ^\circ \Biggr) ,
$$
$$
R_{a,b,(p'-2)/2}(I_2)|0\rangle ''=c_{(p'-2)/2}\Biggl(
a\prod _{j=1}^{p'-1} (ab+\epsilon [j][j+1])^{1/2} \, |p'-1\rangle ^\circ
+(ab+\epsilon [2])^{1/2}|1\rangle ^\circ  $$
$$
+{\rm i} (ab+\epsilon [2])^{1/2} |p'-1\rangle ^\circ +{\rm i}
(ab+\epsilon [p'-2][p'-1])^{1/2}\, |p'-3 \rangle ^\circ \Biggr) .
$$
If
$$
a\prod _{j=1}^{p'-1} (ab+\epsilon [j][j+1])^{1/2}+{\rm i}
(ab+\epsilon [2])^{1/2}=0, \eqno (56) $$
$$
(ab+\epsilon [2])^{1/2}\prod _{j=1}^{p'-1} (ab+\epsilon [j][j+1])^{1/2}
={\rm i}b, \eqno (57)
$$
then
$$
R_{a,b,(p'-2)/2}(I_2)|p'-1\rangle =\frac {-bc_{p'/2}}{\prod _{j=1}^{p'-1}
(ab+\epsilon [j][j+1])^{1/2}}\, |0\rangle ' , $$
$$
R_{a,b,(p'-2)/2}(I_2)|0\rangle '=\frac {{\rm i}(ab+\epsilon [2])^{1/2}}
{q^{(p'-2)/2}-q^{-(p'-2)/2}}\, |1\rangle '+c|p'-1\rangle ', $$
$$
R_{a,b,(p'-2)/2}(I_2)|0\rangle ''=\frac {{\rm i}(ab+\epsilon [2])^{1/2}}
{q^{(p'-2)/2}-q^{-(p'-2)/2}}\, |1\rangle '',
$$
where $c$ is a nonzero coefficient easily determined from the above formulas.
Hence, the subspaces ${\cal H}_1+{\bf C}|p'-1\rangle$ and ${\cal H}_2$
of the representation space are invariant with respect to the representation
$R_{a,b,(p'-2)/2}$ (we denote these subrepresentations by
$R^1_{a,b,(p'-2)/2}$ and $R^2_{a,b,(p'-2)/2}$, respectively). Remark that
$$
{\rm dim}\ {\cal H}_1+{\bf C}|p'-1\rangle ={\textstyle \frac 12}(p'+1) ,
\ \ \ {\rm dim}\ {\cal H}_2 ={\textstyle \frac 12}(p'-1) .
$$
If
$$
a\prod _{j=1}^{p'-1} (ab+\epsilon [j][j+1])^{1/2}-{\rm i}
(ab+\epsilon [2])^{1/2}=0, \eqno (58) $$
$$
(ab+\epsilon [2])^{1/2}\prod _{j=1}^{p'-1} (ab+\epsilon [j][j+1])^{1/2}
=-{\rm i}b, \eqno (59)
$$
then
$$
R_{a,b,(p'-2)/2}(I_2)|p'-1\rangle =\frac {-bc_{p'/2}}{\prod _{j=1}^{p'-1}
(ab+\epsilon [j][j+1])^{1/2}}\, |0\rangle '' , $$
$$
R_{a,b,(p'-2)/2}(I_2)|0\rangle '=\frac {{\rm i}(ab+\epsilon [2])^{1/2}}
{q^{(p'-2)/2}-q^{-(p'-2)/2}}\, |1\rangle ', $$
$$
R_{a,b,(p'-2)/2}(I_2)|0\rangle ''=\frac {{\rm i}(ab+\epsilon [2])^{1/2}}
{q^{(p'-2)/2}-q^{-(p'-2)/2}}\, |1\rangle ''+c|p'-1\rangle ,
$$
where $c$ is a nonzero coefficient. Hence, now
the subspaces ${\cal H}_1$ and ${\cal H}_2+{\bf C}|p'-1\rangle$
of the representation space are invariant. We denote the
subrepresentations on these subspaces by ${\hat R}^1_{a,b,(p'-2)/2}$ and
${\hat R}^2_{a,b,(p'-2)/2}$, respectively). Note that the representation
${\hat R}^1_{a,b,(p'-2)/2}$ is not equivalent to $R^2_{a,b,(p'-2)/2}$
(and the representation
${\hat R}^2_{a,b,(p'-2)/2}$ is not equivalent to $R^1_{a,b,(p'-2)/2}$)
since the parameters $a$ and $b$ determining these representations
satisfy different equations.

If $a$ and $b$ do not satisfy the relations (56) and (57) or the relations
(58) and (59), then the representation $R_{a,b,(p'-2)/2}$ is irreducible.

Let now $p'$ be odd and $\lambda =-q^{(p'-2)/2}$. In this case,
the representation $R_{a,b,-(p'-2)/2}$ is irreducible if
$a$ and $b$ do not satisfy the relations (56) and (57) or the relations
(58) and (59). If
$a$ and $b$ satisfy the relations (56) and (57), then
$R_{a,b,-(p'-2)/2}$ is a reducible representation
and decomposes into the direct sum of two subrepresentations
acting on the subspaces ${\cal H}_1+{\bf C}|p'-1\rangle$ and ${\cal H}_2$.
These subrepresentations are denoted by $R^1_{a,b,-(p'-2)/2}$ and
$R^2_{a,b,-(p'-2)/2}$, respectively, and are determined as
$$
R^i_{a,b,-(p'-2)/2}(I_1)=-R^i_{a,b,(p'-2)/2}(I_1),\ \ \
R^i_{a,b,-(p'-2)/2}(I_2)=-R^i_{a,b,(p'-2)/2}(I_2),\ \ \ i=1,2.
$$
Similarly, if $a$ and $b$ satisfy the relations (58) and (59), then
$R_{a,b,-(p'-2)/2}$ is a reducible representation
and decomposes into the direct sum of two subrepresentations
acting on the subspaces ${\cal H}_1$ and ${\cal H}_2+{\bf C}|p'-1\rangle$.
These subrepresentations are denoted by ${\hat R}^1_{a,b,-(p'-2)/2}$ and
${\hat R}^2_{a,b,-(p'-2)/2}$, respectively, and are determined as
$$
{\hat R}^i_{a,b,-(p'-2)/2}(I_1)=-{\hat R}^i_{a,b,(p'-2)/2}(I_1),\ \ \
{\hat R}^i_{a,b,-(p'-2)/2}(I_2)=-{\hat R}^i_{a,b,(p'-2)/2}(I_2),\ \ \ i=1,2.
$$

{\bf Proposition 11.} {\it Let the conditions (56) and (57) are satisfied.
Then the representations $R^1_{a,b,(p'-2)/2}$, $R^2_{a,b,(p'-2)/2}$,
$R^1_{a,b,-(p'-2)/2}$ and $R^2_{a,b,-(p'-2)/2}$
are irreducible and pairwise nonequivalent. If the conditions (58) and (59)
are satisfied, then
the representations ${\hat R}^1_{a,b,(p'-2)/2}$, ${\hat R}^2_{a,b,(p'-2)/2}$,
${\hat R}^1_{a,b,-(p'-2)/2}$ and ${\hat R}^2_{a,b,-(p'-2)/2}$
are irreducible and pairwise nonequivalent.}
\medskip

{\sl Proof} is similar to that of the previous propositions and we omit it.
\bigskip

\noindent
{\sf X. OTHER REPRESENTATIONS OF
$U_q({\rm so}_3)$ FOR $q$ A ROOT OF UNITY}
\medskip

In the previous section we described irreducible representations of
$U_q({\rm so}_3)$ obtained from irreducible representations of the algebra
${\hat U}_q({\rm sl}_2)$ for $q$ a root of unity. However, at $q$ a root of
unity the algebra $U_q({\rm so}_3)$ has irreducible representations which
cannot be derived from those of ${\hat U}_q({\rm sl}_2)$. They are obtained
as irreducible components of the representations $Q_\lambda$ from Section VII
when one put $q$ equal to a root of unity. We describe these representations
of $U_q({\rm so}_3)$ in this section.

Let $\lambda =q^\tau$ be a nonzero complex number such that $0\le {\rm Re}\,
\tau <1$ and let ${\cal H}$ be the $p'$-dimensional complex vector space
with basis
$$
| m\rangle ,\ \ \ m=0,1,2,\cdots ,p'-1.
$$
We define on this space the operators $Q'_\lambda (I_1)$ and
$Q'_\lambda (I_2)$ determined by the formulas
$$
Q'_\lambda (I_1) | m\rangle =\frac {\lambda q^m+\lambda ^{-1}q^{-m}}
{q-q^{-1}}\, | m\rangle , $$
$$
Q'_\lambda (I_2) | 0\rangle =\frac 1{q-q^{-1}}\, |1\rangle +
\frac 1{q-q^{-1}}\, |p'-1\rangle ,$$
$$
Q'_\lambda (I_2) | p'-1\rangle =\frac 1{q-q^{-1}}\, |p'-2\rangle +
\frac 1{q-q^{-1}}\, |0\rangle ,$$
$$
Q'_\lambda (I_2) | m\rangle =\frac 1{q-q^{-1}}\, |m-1\rangle +
\frac 1{q-q^{-1}}\, |m+1\rangle , \ \ \ m\ne 0,p'-1.
$$
A direct computation shows that these operators satisfy the relations (7)
and (8) and hence determine a representation of $U_q({\rm so}_3)$
which will be denoted by $Q'_\lambda$.
\medskip

{\bf Theorem 6.} {\it If $\lambda \ne 1$ and $\lambda \ne q^{1/2}$, then
the representation $Q'_\lambda$ is irreducible.}
\medskip

{\sl Proof} of this proposition is the same as that of the first
part of Proposition 3.
\medskip

The representations $Q'_1$ and $Q'_{\sqrt q}$ are studied in the same
way as the representations $Q_1$ and $Q_{\sqrt q}$ in Section VII. This study
leads to the irreducible representations of $U_q({\rm so}_3)$ which
are described below. (Note that the description of these representations
for $p'$ even and for $p'$ odd is deferent.)

Let $p'$ be odd. We denote by ${\cal H}_r$ and ${\cal H}_s$, $r=\frac 12
(p'+1)$, $s=\frac 12 (p'-1)$, the complex vector spaces with the bases
$$
|0\rangle ,\ |1\rangle ,\ |2\rangle ,\ \cdots \ ,\
|{\textstyle \frac 12}(p'-1)\rangle
$$
and
$$
|1\rangle ,\ |2\rangle ,\ \cdots \ ,\ |{\textstyle \frac 12}(p'-1)\rangle ,
$$
respectively. Four representations $Q_1^{\pm ,\pm}$ act on the space
${\cal H}_r$ and are given by the formulas
$$
Q_1^{+,\pm}(I_1)|m\rangle =\frac {q^m+q^{-m}}{q-q^{-1}} \,
|m\rangle ,\ \ \ m=0,1,2,\cdots , {\textstyle \frac 12}(p'-1), \eqno (60)
$$
$$
Q_1^{+,\pm }(I_2)| {\textstyle \frac 12}(p'-1)\rangle =\pm \frac 1{q-q^{-1}}
\, |{\textstyle \frac 12} (p'-1)\rangle
+ \frac 1{q-q^{-1}} \,
|{\textstyle \frac 12} (p'-3)\rangle , \eqno (61) $$
$$
Q_1^{+,\pm }(I_2)| m\rangle = \frac 1{q-q^{-1}} \,
|m+1\rangle + \frac 1{q-q^{-1}} \, |m-1\rangle , \ \ \
m<{\textstyle \frac 12} (p'-1),
\eqno (62)
$$
and by the formulas
$$
Q_1^{-,\pm}(I_1)|m\rangle =- \frac {q^m+q^{-m}}{q-q^{-1}} \,
|m\rangle ,\ \ \ m=0,1,2,\cdots , {\textstyle \frac 12}(p'-1), \eqno (63)
$$
$$
Q_1^{-,\pm }(I_2):=Q_1^{+,\pm }(I_2). \eqno (64)
$$
Note that the upper sign corresponds to the representations $Q_1^{+,+}$ and
$Q_1^{-,+}$ and the lower sign to the representations $Q_1^{+,-}$ and
$Q_1^{-,-}$.

On the space ${\cal H}_s$, four representations
${\hat Q}_1^{\pm ,\pm}$ act by the corresponding formulas (60)--(64), but
now $m$ runs over the values $1, 2, 3,\cdots ,\frac 12 (p'-1)$.

Let now ${\cal H}'_r$ and ${\cal H}'_s$, $r=\frac 12
(p'+1)$, $s=\frac 12 (p'-1)$, be the complex vector spaces with the bases
$$
| m+{\textstyle \frac 12}\rangle ,\ \ \ m=0,1,2,\cdots ,
{\textstyle \frac 12}(p'-1) ,
$$
and
$$
| m+{\textstyle \frac 12}\rangle ,\ \ \ m=0,1,2,\cdots ,
{\textstyle \frac 12}(p'-3) ,
$$
respectively. The four representations $Q_{\sqrt q}^{\pm ,\pm}$ act on the
space ${\cal H}'_r$ and are given by the formulas
$$
Q_{\sqrt q}^{+ ,\pm}(I_1) |m+{\textstyle \frac 12}\rangle =
\frac {q^{m+1/2}+q^{-m-1/2}}{q-q^{-1}}\, |m+{\textstyle \frac 12}\rangle ,\ \
m=0,1,2,\cdots , {\textstyle \frac 12}(p'-1) , \eqno (65)  $$
$$
Q_{\sqrt q}^{+ ,\pm}(I_2) |{\textstyle \frac 12}\rangle =\pm
\frac 1{q-q^{-1}}\, |{\textstyle \frac 12}\rangle
+\frac 1{q-q^{-1}}\, |{\textstyle \frac 32}\rangle  , \eqno (66) $$
$$
Q_{\sqrt q}^{+ ,\pm}(I_2) |m+{\textstyle \frac 12}\rangle =\frac 1{q-q^{-1}}
|m+{\textstyle \frac 32}\rangle +\frac 1{q-q^{-1}}
|m-{\textstyle \frac 12}\rangle ,\ \ \
m\ne 0,  \eqno (67)
$$
where $|m+{\textstyle \frac 32}\rangle \equiv 0$ if $m=\frac 12 (p'-1)$,
and by the formulas
$$
Q_{\sqrt q}^{- ,\pm}(I_1) |m+{\textstyle \frac 12}\rangle =-
\frac {q^{m+1/2}+q^{-m-1/2}}{q-q^{-1}}\, |m+{\textstyle \frac 12}\rangle ,\ \
m=0,1,2,\cdots , {\textstyle \frac 12}(p'-1) , \eqno (68)  $$
$$
Q_{\sqrt q}^{- ,\pm}(I_2) :=Q_{\sqrt q}^{+ ,\pm}(I_2). \eqno (69)
$$
On the space ${\cal H}'_s$, four representations
${\breve Q}_{\sqrt q}^{\pm ,\pm}$ act by the corresponding formulas
(65)--(69),
but now $m$ runs through the values $0,1,2,\cdots ,\frac 12 (p'-3)$.

Let now $p'$ be even. We denote by ${\cal H}_r$ and ${\cal H}_s$,
$r=\frac 12 (p'+2)$, $s=\frac 12 (p'-2)$, the complex vector spaces with
the bases
$$
| 0\rangle ,\ |1\rangle ,\ |2\rangle ,\ \cdots \ ,\
|{\textstyle \frac 12} p'\rangle
$$
and
$$
 |1\rangle ,\ |2\rangle ,\ \cdots \ ,\
|{\textstyle \frac 12} (p'-2)\rangle ,
$$
respectively. The representations $Q_1^{1,\pm}$ and
$Q_1^{2,\pm}$ act on ${\cal H}_r$ and ${\cal H}_s$,
respectively, which are given by the formulas
$$
Q_1^{i,\pm} (I_1)| m\rangle =\pm \frac {q^m+q^{-m}}{q-q^{-1}}\,
| m\rangle ,\ \ \ i=1,2 , $$
$$
Q_1^{i,\pm} (I_2)| m\rangle =\frac 1{q-q^{-1}}\, | m+1\rangle +
\frac 1{q-q^{-1}}\, | m-1\rangle , \ \ \ i=1,2 ,
$$
where $|m+1\rangle $ or $|m-1\rangle$ must be put equal to 0 if the
corresponding vector does not exist.

Let ${\cal H}_{p'/2}$ be the complex vector space with the basis
$$
|m+{\textstyle \frac 12} \rangle ,\ \ \ m=0,1,2,\cdots ,
{\textstyle \frac 12} (p'-2).
$$
Four representations ${\hat Q}_{\sqrt q}^{\pm ,\pm}$ act on this space
which are given by the formulas
$$
{\hat Q}_{\sqrt q}^{+ ,\pm}(I_1)|m+{\textstyle \frac 12} \rangle =
\frac {q^{m+1/2}+q^{-m-1/2}}{q-q^{-1}}| m+{\textstyle \frac 12}\rangle ,$$
$$
{\hat Q}_{\sqrt q}^{+ ,\pm}(I_2)|{\textstyle \frac 12} \rangle =\pm
\frac 1{q-q^{-1}}| {\textstyle \frac 12}\rangle
+\frac 1{q-q^{-1}}| {\textstyle \frac 32}\rangle , $$
$$
{\hat Q}_{\sqrt q}^{+ ,\pm}(I_2)|{\textstyle \frac 12}(p'-2) \rangle =\pm
\frac 1{q-q^{-1}}| {\textstyle \frac 12}(p'-2)\rangle
+\frac 1{q-q^{-1}}| {\textstyle \frac 12}(p'-4)\rangle , $$
$$
{\hat Q}_{\sqrt q}^{+ ,\pm}(I_2)|{\textstyle m+\frac 12} \rangle =
\frac 1{q-q^{-1}}| {\textstyle m-\frac 12}\rangle
+\frac 1{q-q^{-1}}| {\textstyle m+\frac 32}\rangle ,\ \ m\ne
{\textstyle \frac 12}, {\textstyle \frac 12}(p'-2) ,
$$
and by the formulas
$$
{\hat Q}_{\sqrt q}^{- ,\pm}(I_1)|m+{\textstyle \frac 12} \rangle =-
\frac {q^{m+1/2}+q^{-m-1/2}}{q-q^{-1}}| m+{\textstyle \frac 12}\rangle ,$$
$$
{\hat Q}_{\sqrt q}^{- ,\pm}(I_2) ={\hat Q}_{\sqrt q}^{+ ,\pm}(I_2).
$$

Let us mention peculiarities of the representations described above.
The operators $Q_1^{\pm ,\pm}(I_2)$, $Q_{\sqrt q}^{\pm ,\pm}(I_2)$,
${\breve Q}_{\sqrt q}^{\pm ,\pm}$ and
${\hat Q}_{\sqrt q}^{\pm ,\pm}(I_2)$ have nonzero diagonal matrix
elements and nonzero traces. Moreover, the operators
${\hat Q}_{\sqrt q}^{\pm ,\pm}(I_2)$ have two such diagonal elements.
Spectra of the operators
$Q_1^{\pm ,\pm}(I_1)$, $Q_{\sqrt q}^{\pm ,\pm}(I_1)$,
$Q_1^{1 ,\pm}(I_1)$, $Q_1^{2 ,\pm}(I_1)$ and
${\hat Q}_{\sqrt q}^{\pm ,\pm}(I_1)$ are not symmetric with respect to
the zero point.
\medskip

{\bf Proposition 12.} {\it The representations
$Q_1^{\pm ,\pm}$, $Q_{\sqrt q}^{\pm ,\pm}$,
${\breve Q}_{\sqrt q}^{\pm ,\pm}$, $Q_1^{1 ,\pm}$, $Q_1^{2 ,\pm}$,
${\hat Q}_{\sqrt q}^{\pm ,\pm}$ are irreducible and pairwise
nonequivalent. No representation $Q'_\lambda$ is equivalent to any of
these representations.}
\medskip

{\sl Proof} is the same as that of Proposition 3.
\bigskip

\noindent
{\sf ACKNOWLEDGMENTS}
\medskip

The second co-author (A.U.K.) was supported in part by Civilian Reaserch and
Development Foundation Grant No. UP1-309. This co-author would like to
thank Department of Mathematics of FNSPE, Czech Technical University,
for the hospitality.
\bigskip
\bigskip

\noindent
1. R. M. Santilli, {\it Nouvo Cim.} {\bf 51}, 570 (1967).

\noindent
2. R. M. Santilli, {\it Hadronic J.} {\bf 1}, 574 (1978).

\noindent
3. R. M. Santilli, {\it Hadronic J.} {\bf 4}, 1166 (1981).

\noindent
4. D. B. Fairlie, {\it J. Phys. A} {\bf 23}, L183 (1990).

\noindent
5. M. Odesski, Funct. Anal. Appl. {\bf 20}, No. 2, 78 (1986).

\noindent
6. Yu. S. Samoilenko and L. B. Turovska, {\it Preprint}, Kiev Institute
of Mathematics, 1996 ({\it Reps. Math. Phys.}, in press).

\noindent
7. O. V. Bagro and S. A. Kruglyak, {\it Preprint}, Kiev, 1996.

\noindent
8. M. Havl\'{i}\v cek, A. Klimyk and E. Pelantov\'a,
Czech. J. Phys. {\bf 47}, 13 (1997).

\noindent
9. M. Noumi, {\it Adv. Math.} {\bf 123}, 11 (1996).

\noindent
10. J. Dixmier, {\it Alg\'ebres Enveloppantes}, Gauthier--Villars, Paris, 1974.

\noindent
11. A. Klimyk and K. Schm\"udgen, {\it Quantum Groups and Their
Representations}, Springer, Berlin, 1997.

\noindent
12. I. M. Burban and A. U. Klimyk, {\it J. Phys. A} {\bf 26}, 2139 (1993).

\noindent
13. A. M. Gavrilik and A. U. Klimyk, {\it J. Math. Phys.}
{\bf 35}, 3670 (1994).

\end{document}